\documentclass[%
superscriptaddress,
groupedaddress,
 amsmath,amssymb,
 aps,
 onecolumn,
physrev,
floatfix,
]{revtex4-2}

\usepackage{graphicx}%
\usepackage{amsmath}
\usepackage{amssymb}
\usepackage{amsfonts}
\usepackage[style=base]{subcaption}
\usepackage{hyperref}%
\usepackage{placeins}

\usepackage{bm}

\newcommand{\ABS}[1]{\ensuremath{\lvert {#1} \rvert}}
\newcommand{\dpone}[2]{\ensuremath{\displaystyle\frac{\partial {#1}}{\partial {#2}}}}

\newcommand{\bD}{\ensuremath{\bm{D}}}

\newcommand{\bG}{\ensuremath{\bm{G}}}

\newcommand{\bI}{\ensuremath{\bm{I}}}

\newcommand{\bK}{\ensuremath{\bm{K}}}
\newcommand{\bL}{\ensuremath{\bm{L}}}

\newcommand{\bQ}{\ensuremath{\bm{Q}}}

\newcommand{\bS}{\ensuremath{\bm{S}}}
\newcommand{\bT}{\ensuremath{\bm{T}}}

\newcommand{\bW}{\ensuremath{\bm{W}}}

\newcommand{\bd}{\ensuremath{\bm{d}}}
\newcommand{\be}{\ensuremath{\bm{e}}}
\newcommand{\bff}{\ensuremath{\bm{f}}}
\newcommand{\bg}{\ensuremath{\bm{g}}}

\newcommand{\bn}{\ensuremath{\bm{n}}}

\newcommand{\br}{\ensuremath{\bm{r}}}
\newcommand{\bs}{\ensuremath{\bm{s}}}
\newcommand{\bt}{\ensuremath{\bm{t}}}
\newcommand{\bu}{\ensuremath{\bm{u}}}

\newcommand{\bx}{\ensuremath{\bm{x}}}
\newcommand{\by}{\ensuremath{\bm{y}}}

\newcommand{\bbR}{\ensuremath{\mathbb{R}}}
\newcommand{\bbRv}{\ensuremath{\mathbb{R}^3}}
\newcommand{\bbRt}{\ensuremath{\mathbb{R}^{3\times 3}}}
\newcommand{\Tcal}{\ensuremath{\mathcal{T}}}

\newcommand{\SM}{\ensuremath{\text{S}\to\text{M}}}
\newcommand{\SL}{\ensuremath{\text{S}\to\text{L}}}
\newcommand{\ST}{\ensuremath{\text{S}\to\text{T}}}
\newcommand{\MM}{\ensuremath{\text{M}\to\text{M}}}
\newcommand{\ML}{\ensuremath{\text{M}\to\text{L}}}
\newcommand{\MT}{\ensuremath{\text{M}\to\text{T}}}

\newcommand{\LL}{\ensuremath{\text{L}\to\text{L}}}
\newcommand{\LT}{\ensuremath{\text{L}\to\text{T}}}
\newcommand{\MLtoML}{\ensuremath{\text{ML-tree}}}

\newcommand{\bsigma}{\ensuremath{\bm{\sigma}}}
\newcommand{\bomega}{\ensuremath{\bm{\omega}}}
\newcommand{\bOmega}{\ensuremath{\bm{\Omega}}}

\begin{document}

\title{Kernel Aggregated Fast Multipole Method: Efficient summation of Laplace and Stokes kernel functions}

\author{Wen Yan}
\email{wyan@flatironinstitute.org}
\author{Robert Blackwell}
\address{Center for Computational Biology, Flatiron Institute, Simons Foundation, New York, NY 10010}

\date{\today}%

\begin{abstract}
  Many different simulation methods for Stokes flow problems involve a common computationally intense task---the summation of a kernel function over \(O(N^2)\) pairs of points.
  One popular technique is the Kernel Independent Fast Multipole Method (KIFMM), which constructs a spatial adaptive octree for all points and places a small number of equivalent multipole and local equivalent points around each octree box, and completes the kernel sum with \(O(N)\) cost, using these equivalent points.
  Simpler kernels can be used between these equivalent points to improve the efficiency of KIFMM.
  Here we present further extensions and applications to this idea, to enable efficient summations and flexible boundary conditions for various kernels.
  We call our method the Kernel Aggregated Fast Multipole Method (KAFMM), because it uses different kernel functions at different stages of octree traversal.
  We have implemented our method as an open-source software library \texttt{STKFMM} based on the high performance library \texttt{PVFMM}, with support for Laplace kernels, the Stokeslet, regularized Stokeslet, Rotne-Prager-Yamakawa (RPY) tensor, and the Stokes double-layer and traction operators.
  Open and periodic boundary conditions are supported for all kernels, and the no-slip wall boundary condition is supported for the Stokeslet and RPY tensor.
  The package is designed to be ready-to-use as well as being readily extensible to additional kernels.
  % \keywords{Fast Multipole Method, Stokeslet, Periodic Boundary Conditions}
\end{abstract}

\maketitle

\section{\label{sec:intro}Introduction}
Solving the Stokes flow field in the presence of multiphase and complex geometry is a classic problem appearing in many interesting areas, including active swimmers, colloidal dynamics, and rheology of complex fluids.
Many numerical methods have been developed, such as Stokesian Dynamics (SD), Force Coupling Method (FCM), and Boundary Integral (BI).
Each numerical method has certain advantages and disadvantages compared to one another, but they all rely on the classic kernel summation problem.
This problem is defined by a set of sources (S) with locations \(\by^\beta\), and a kernel function \( \bK \), from which certain quantities \(\bg^\alpha\) are to be calculated at target points (T) located at \(\bx^\alpha\):
\begin{align}\label{eq:kernelsumlinear}
  g_i^\alpha  =\sum_{\beta\leq N_S} \bK_{ij}\left(\bx^\alpha,\by^\beta\right)s_j^\beta, \quad \forall \alpha\leq N_T.
\end{align}
Here \(N_S\) is the number of source points, \(N_T\) is the number of target points, and
\(\bs^\beta\) denotes the (vectorial or scalar) source values. The \(\bg^\alpha\) values also may be either vectorial or scalar.
For clarity, we will always use \(\bx\) for a target point and \(\by\) for a source point, although in many applications they may actually be the same set of points.

The main challenge of Eq.~(\ref{eq:kernelsumlinear}) is its quadratic \(O(N_S N_T)\) computational cost, since every source-target pair must naively be evaluated.
This appears necessary because the kernel function \(\bK\) often decays slowly and cannot be truncated with a cut-off radius.
This is the case when \(\bK\) is the Green's function of the Laplace or Stokes equations.
To reduce the computational complexity of the prohibitive \(O(N_S N_T)\) summation, one popular technique is the use of Ewald-type methods.
In Ewald-type methods, the kernel function \(\bK\) is split into a fast-decaying short-range part and a slowly varying long-range part via a splitting parameter \(\eta\).
By properly setting \(\eta\), the short-range part can be truncated with a very short cut-off radius, while the long-range part is efficiently summed on a spatial regular mesh using a Fast Fourier Transform (FFT).
In this way, the overall complexity is reduced to \(O(N_S+N_T)+O(M\log M)\), where \(M\) is the number of regular mesh grid points used in the FFT operation.
Methods in this fashion include Particle-Mesh-Ewald (PME), Spectral Ewald, Accelerated Stokesian Dynamics, and others \cite{SwanBrady2011Anisotropicdiffusionconfined,WangBrady2016SpectralEwaldAcceleration,YanBrady2016behavioractivediffusiophoretic}.
While these methods were initially applicable only to periodic boundary conditions, they have recently been extended to open boundary conditions and partially periodic boundary conditions
\cite{Tornberg2015Ewaldsumssingly,KlintebergShamshirgarEtAl2017FastEwaldsummation,LindboTornberg2012Fastspectrallyaccurate,ShamshirgarTornberg2017SpectralEwaldmethod,SrinivasanTornberg2018FastEwaldSummation}.

Another choice is the fast multipole method (FMM).
The original FMM works by expanding the kernel into a truncated series of spherical harmonic functions and grouping source and target points into a spatial adaptive octree grid.\cite{GreengardRokhlin1987fastalgorithmparticle}
The interactions between source and target pairs in non-adjacent octree boxes are then represented by the spherical harmonic functions instead of being evaluated pairwise.
In this way, the overall computational cost is reduced to \(O(N_S+N_T)\), since the maximum number of interaction pairs for each point is limited by the number of adjacent octree boxes and the order of the spherical harmonic functions.
FMM also typically has an advantage over FFT-based methods when handling highly nonuniform distributions of points thanks to the adaptive octree.

The kernel independent fast multipole method (KIFMM) is a variant of FMM which replaces the spherical harmonic functions with a set of equivalent points surrounding each octree box.
The far-field generated by source points in an octree leaf is approximated as if generated by the multipole equivalent points surrounding this box (M points, also called ``upward points'').
Meanwhile, the field induced by faraway source points ``felt'' by target points in a leaf box is approximated by a set of local (L, also called ``downward'') equivalent points.
Each multipole or local equivalent point generates a field using the same kernel function \(\bK\).
This method is called ``kernel independent'' because the spherical harmonic expansion for different kernel functions is eliminated and the same code can be used for different kernel functions \(\bK\).
More details about the KIFMM algorithm can be found in \cite{YingBirosEtAl2004kernelindependentadaptive}.
Although KIFMM was initially developed for open boundary conditions, it since has been extended to allow for fully or partially periodic boundary conditions \cite{YanShelley2018Flexiblyimposingperiodicity}.

Despite its linear theoretical computational cost and spatial adaptivity, there are a few cases where the advantage of KIFMM has not been widely recognized.
For instance, computations of the flow field generated by a spherical particle in a Stokes fluid often use the Rotne-Prager-Yamakawa (RPY) kernel \(\bu=\left(1+\frac{1}{6}b^2\nabla^2\right)\bG\bff\), where \( b \) is the particle radius and \(\bG\) is the Stokeslet.
In some other cases, the regularized Stokeslet $G_{ij}^\epsilon = \frac{1}{8\pi}\left[\frac{r^{2}+2 \epsilon^{2}}{\left(r^{2}+\epsilon^{2}\right)^{3 / 2}} \delta_{i j}+\frac{1}{\left(r^{2}+\epsilon^{2}\right)^{3 / 2}} r_i r_j \right]$ is preferred, where $\epsilon$ is the source spreading length scale.
In both cases the kernels depend linearly on the force \(\bff\) but nonlinearly on the particle radius \(b\) or the regularization parameter $\epsilon$.
This nonlinearity sometimes poses a difficulty for KIFMM, especially when $b$ and $\epsilon$ are different for each source point.
The RPY kernel has been discussed in the past by conversion to the combination of a few harmonic FMMs \cite{LiangGimbutasEtAl2013FastMultipoleMethod}, by a direct expansion into spherical harmonic functions \cite{GuanChengEtAl2018RPYFMMParallelAdaptive}, or by barycentric Lagrange far-field expansions\cite{wangKernelIndependentTreecodeGeneral2021},
but all these methods require sometimes involved derivations or expansions of the RPY kernel function.
The regularized Stokeslet has been directly applied in KIFMM\cite{rostamiKernelindependentFastMultipole2016}, but polydispersity (different $\epsilon$ for different points) and accuracy better than $10^{-5}$ were not demonstrated.

In this work we extend the core idea of `equivalent sources' in KIFMM, to handle various kernels that are more complicated than the Stokeslet, including the RPY and regularized variants.
In KIFMM, sometimes simpler sources can be placed on M and L points to improve the efficiency of the tree traversal stage.
For example, charges, instead of dipoles, can be placed on M and L points when evaluating the potential generated by Laplace dipoles.
This idea has been utilized by \texttt{PVFMM} \cite{MalhotraBiros2015PVFMMParallelKernel} and some other implementations, but limited to cases such as Laplace monopoles replacing dipoles and Stokes single layer sources replacing double layer sources.
We extend this idea to various other kernels that are important in applications.
The cost of tree traversal can be minimized and, more importantly, nonlinearity in regularized Stokeslet and RPY kernels can be conveniently handled.
In addition, full or partial periodic boundary conditions can be conveniently imposed using methods developed in our previous work \cite{YanShelley2018Flexiblyimposingperiodicity,YanShelley2018UniversalImageSystems}, because we only need the M-to-L translation kernel to impose the boundary condition and different kernel summation problems may share the same M-to-L kernel and use the same precomputed periodization operator.

In Section~\ref{sec:method} we discuss the choices of different kernels for an assortment of kernel summation problems.
In Section~\ref{sec:bench} we demonstrate the convergence of the KAFMM algorithm by reporting a representative selection of numerical convergence and timing results using our ready-to-use and open-source software package, \texttt{STKFMM}. \footnote{https://github.com/wenyan4work/STKFMM}
The package includes several kernel functions useful in solving Laplace and Stokes problems.
The efficient, scalable, and flexible design of our package is enabled by the highly optimized \texttt{PVFMM} library \cite{MalhotraBiros2015PVFMMParallelKernel}.
We include some brief parallel scaling benchmark results in Appendix \ref{app:scaling}.
We conclude with some brief remarks in Section~\ref{sec:conclusion}.

\section{\label{sec:method}Methodology}
We first rewrite the kernel summation problem Eq.~(\ref{eq:kernelsumlinear}) more generally,
allowing the kernel \(\bK\) to be a nonlinear function of the source values \(\bs^\beta\):
\begin{align}\label{eq:kernelsum}
  \bg^\alpha  =\sum_{\beta\leq N_S} \bK\left(\bx^\alpha,\by^\beta, \bs^\beta\right), \quad \forall \alpha\leq N_T.
\end{align}

For a function \(f(\br)\), where \(\br=\bx-\by\), we have $\dpone{f}{x_j} = -\dpone{f}{y_j}$, and $\displaystyle \frac{\partial^2 f}{\partial x_j \partial x_k} = \frac{\partial^2 f}{\partial y_j \partial y_k}$.
Naturally, \( \displaystyle \frac{\partial^2 f}{\partial x_j \partial x_k}\) is a symmetric tensor when \(f(\br)\) is a smooth function.

For every kernel \(\bK\) we define its source dimension \(k_s\) (target dimension \(k_t\)) to be the dimension of its source value \(\bs\) (target value \(\bg\)).
\(\bK\) is therefore a (possibly nonlinear) mapping \(\bbR^{k_s}\to\bbR^{k_t}\) and we say \(\bK\) has dimensions \(k_s\times k_t\).
For example, a scalar kernel \(K\) has dimensions \(1\times 1\), the gradient \(\nabla K\) has dimensions \(1\times 3\) in 3D space, and the gradient of gradient \(\nabla\nabla K\) has dimensions \(1\times 6\) due to the symmetry of \(\nabla\nabla K\).

In many cases, we want to know not only the target value \(\bg\), but also the gradient \(\nabla\bg\), the gradient of the gradient \(\nabla\nabla\bg\), or the Laplacian \(\nabla^2\bg\).
We define ``combined'' kernels with the symbol \(\oplus\), for example, \(\bK\oplus\nabla\bK\).
If \(\bK\) has dimensions \(m\times n\), then \(\nabla\bK\) has dimensions \(m\times 3n\) in 3D space,
and the combined kernel \(\bK\oplus\nabla\bK\) has dimensions \(m \times (n+3n)\).
For the Stokeslet \(\bG\) with dimensions \(3\times 3\), for instance, the combined kernel \(\bG\oplus\nabla^2\bG\) has dimensions \(3\times 6\) and maps the force at a source point to the velocity and Laplacian of velocity at a target point.

Further, we shall write the kernel summation problem Eq.~(\ref{eq:kernelsum}) in a more compact form as \(\bs \to \left[\bg\right]\).
For the combined kernel \(\bK\oplus\nabla\bK\), for example, we write
\begin{align}\label{eq:defcombsum}
  \bs \to \left[\bg,\nabla\bg\right],
\end{align}
which refers to the problem where both \(\bg\) and its gradient \(\nabla\bg\) are calculated at each target point:
\begin{align}
  \bg^\alpha       & =\sum_{\beta\leq N_S} \bK\left(\bx^\alpha,\by^\beta, \bs^\beta\right), \quad \forall \alpha\leq N_T,       \\
  \nabla\bg^\alpha & =\sum_{\beta\leq N_S} \nabla\bK\left(\bx^\alpha,\by^\beta, \bs^\beta\right), \quad \forall \alpha\leq N_T.
\end{align}

\begin{figure}[htbp]
  \centering
  \includegraphics[width=\linewidth]{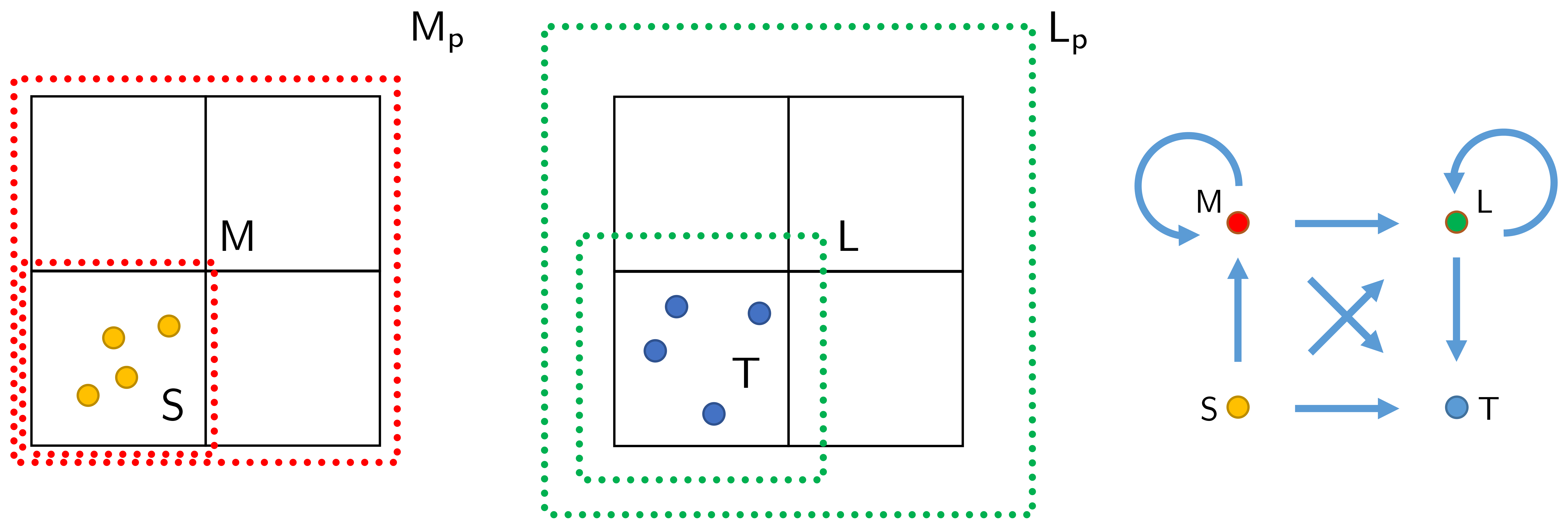}
  \caption{\label{fig:kifmmgraph}
    The KIFMM algorithm graph.
    Left: Some source points (yellow) located in an octree leaf box, which is surrounded by multipole points \(\text{M}\) (red). Further multipole points \(\text{M}_p\) surround the parent box. Middle: Some target points (blue) located in another octree leaf box, which is surrounded by local points \(\text{L}\) (green), with further local points \(\text{L}_p\) corresponding to its parent box.
    The root octree box is not shown.
    The source and target leaf boxes need not belong to the same level in the octree, and are shown with the same size only for clarity.
    Right: Graph of all the possible routes between source and target points. For specific source and target points, the route depends on whether their leaf boxes are adjacent and whether they are at the same level of the adaptive octree.}
\end{figure}

The octree traversal algorithm in KIFMM consists of an upward pass followed by a downward pass.
The upward pass first determines the sources at M points surrounding each leaf box that approximate the field induced faraway by the S points in the leaf box (see Fig.~\ref{fig:kifmmgraph}).
This step is called \SM{}.
This is done by placing a certain number of check points on a surface surrounding this octree box, and than ensure that the sources within this octree box and the equivalent sources on the M points generate the same field on those check points.
Then, the upward pass repeats this process and recursively moves to higher level octree boxes and for each one finds the sources at M points to approximate the field induced by enclosed M points one level lower down.
This step is called \MM{}.
Once the root octree box is reached, the upward pass finishes and the downward pass begins.
In the downward pass, sources on L points are determined similarly as the M points in the upward pass, by placing a certain number of check points on a surface surrounding this octree box, and then ensure that the sources far from this octree box and the equivalent sources on the L points generate the same field on those check points.
Similar to the upward pass, the downward pass recursively traverses all octree boxes down to the lowest leaf level, i.e., \LL{}.
The strengths at L points approximate the field induced inside the associated box by distant sources outside of it---which may include S points, M points, or other L points, depending on the stage of the tree traversal and relative box locations and levels.
Determining M and L sources involves solving inverse problems, i.e., solving for source strengths from field values on check points.
In the original KIFMM \cite{YingBirosEtAl2004kernelindependentadaptive}, Tikhonov regularization is used to stabilize the inversion.
This was replaced by backward stable SVD in \texttt{PVFMM} \cite{MalhotraBiros2015PVFMMParallelKernel} to improve the accuracy and performance.

The kernel summation problem includes all S and T pairs.
In the KIFMM algorithm, any given pair of S and T points will be evaluated along one of the routes in the graph on the right of Fig.~\ref{fig:kifmmgraph}.
If S and T are located in adjacent leaf boxes, the equivalent M and L points cannot be used and the path \ST{} is taken.
If S and T are in separated leaf boxes, then the path connecting them will be one of the many other different possible paths, depending on their relative locations and octree levels.
Possible paths include \(\text{S}\to\text{M}\to\text{L}\to\text{T}\), \(\text{S}\to\text{M}\to\text{M}\to\text{L}\to\text{T}\), \(\text{S}\to\text{L}\to\text{T}\), and so on.
There are no paths consisting of \(\text{L}\to\text{M}\) or \(\text{M}\to\text{S}\), because the paths are all unidirectional from S to T.
The 8 arrows in the graph correspond to 8 different sub-operations to be performed during a step of the KIFMM tree traversal, using a suitably chosen kernel.
We summarize these operations in a kernel table, as shown in Table~\ref{tab:ktabdef}.
The decision about which path the evaluation of a specific pair (S,T) takes is based on U,V,W,X lists constructed for each octree leaf box, which classify the other boxes according to their relative locations \cite{YingBirosEtAl2004kernelindependentadaptive,MalhotraBiros2015PVFMMParallelKernel}.
For a leaf box B, the source points within these four lists must be handled by B itself, because these boxes are not far separated from B.
The contribution from more distant boxes are considered by higher level boxes containing B.
For example, U-list refers to the octree boxes adjacent to B such that direct \ST{} must be evaluated, and V-list refers to the set of the children of the neighbors of the parent of B which are not adjacent to B.
Our method does not modify these lists or how they are used, so we omit those technical details here.

\begin{table}[htbp]
  \centering
  \begin{tabular}{c|c|c|c}
    \hline
      & M   & L   & T   \\
    \hline
    S & \SM & \SL & \ST \\
    M & \MM & \ML & \MT \\
    L & --  & \LL & \LT \\
    \hline
  \end{tabular}
  \caption{\label{tab:ktabdef} The general kernel table of the KIFMM method.
    The kernel from L to M is undefined because there is no such operation in the KIFMM algorithm. }
\end{table}

In the original KIFMM, the same kernel is used for all 8 interactions.
One well-known extension\cite{MalhotraBiros2015PVFMMParallelKernel} to this is when the function value $\bK$ and its derivatives $\nabla\bK$ or $\nabla^2\bK$ are both needed on target points.
For example, when evaluating the kernel summation problem \(\bs\to\left[\bg,\nabla^2\bg\right]\),
we can use the \(\bK\) kernel for many of the interactions and we only need to use the combined kernel for interactions where the target values are computed, namely \ST, \MT, and \LT, as shown in Table~\ref{tab:ktab_example}.
\begin{table}[htbp]
  \centering
  \begin{tabular}{c|c|c|c}
    \hline
      & M       & L       & T                        \\
    \hline
    S & \(\bK\) & \(\bK\) & \(\bK\oplus\nabla^2\bK\) \\
    M & \(\bK\) & \(\bK\) & \(\bK\oplus\nabla^2\bK\) \\
    L & --      & \(\bK\) & \(\bK\oplus\nabla^2\bK\) \\
    \hline
  \end{tabular}
  \caption{\label{tab:ktab_example} The KAFMM kernel table for the kernel summation problem \(\bs\to\left[\bg,\nabla^2\bg\right]\). }
\end{table}
This simplification is possible because the KIFMM algorithm is based on the concept of equivalence:
For S and M points corresponding to an octree leaf box, KIFMM checks the equivalence of the field (or the potential) generated by the S points and that generated by the M points at a specific ``check surface'' surrounding the box.
Based on the theory of the method of fundamental solutions (MFS), if the S and M points generate the same potential on this check surface, they generate the same potential for all points outside this check surface.
Therefore, when checking the equivalence there is no need to check the Laplacians---they are guaranteed to match at all locations outside the check surface if the potentials generated by the S and M points match at the check surface.
The same results also apply to other derivatives of \(\bK\), including \(\nabla\bK\), \(\nabla\nabla\bK\), and higher order derivatives.
Similar results apply to the case when L points are being checked.
This flexibility also works when S points include higher order moments.
For instance, if the S points are Laplace quadrupoles, the M and L points can be Laplace charges (i.e., monopoles) because they have no difficulty in matching the potential generated by Laplace quadrupoles.
This idea has been well known and implemented by many KIFMM packages, for example, PVFMM \cite{MalhotraBiros2015PVFMMParallelKernel}.

Our KAFMM extends this idea.
For a given kernel summation problem, the specific choice of aggregated kernels depends on the \ST{} kernel we want to evaluate.
The most crucial ingredient is to choose the best kernels for tree traversal operations, i.e., the \MM{}, \ML{}, and \LL{} kernels, which we will refer to as the \MLtoML{} kernels for brevity.
Those kernels must be linear even for a problem with a nonlinear \ST{} kernel.
The chosen \MLtoML{} kernels must also have low kernel dimensions to improve the octree traversal efficiency.
In the following sections we discuss the choice of kernels for a variety of Laplace and Stokes kernel summation problems.

\subsection{\label{subsec:laplace}Laplace kernels}
The Laplace kernels are straightforward but we still include these well-known results here as an example of our notation and implementation.
The Laplace single-layer kernel \(L\), which is the Green's function for the Laplace equation,
gives the potential \(\phi\) generated by charge \(q\): \(\phi = Lq\). Explicitly,
\begin{align}
  L(\bx,\by)=\frac{1}{4\pi r},
\end{align}
where \(\br=\bx-\by\) and $r=\ABS{\br}$.
The dipole and quadrupole kernels, \(\bL^D\) and \(\bL^Q\), give the potentials generated by a dipole source \(d_j\in \bbRv\) and a quadrupole source \(Q_{jk}\in \bbRt\):
\begin{align}
  \phi & =L_j^D d_j, \quad L_j^D = -\dpone{L}{x_j},                                      \\
  \phi & =L_{jk}^Q Q_{jk}, \quad L_{jk}^Q=\frac{\partial^2 L}{\partial x_j\partial x_k},
\end{align}
where we use the comma notation defined earlier.
The three kernels \(L, \bL^D, \bL^Q\) are all linear operators applied to the source values.

As mentioned in the previous part, a useful case is when the kernel summation problem has both charges and dipoles on S points, i.e., \([q,\bd]\to\left[\phi,\nabla\phi,\nabla\nabla\phi\right]\).
In this case, there is no need to run the algorithm twice and evaluate \(q\to\left[\phi,\nabla\phi,\nabla\nabla\phi\right]\) and \(\bd\to\left[\phi,\nabla\phi,\nabla\nabla\phi\right]\) separately.
The calculation can be done in a single KAFMM pass, using the kernel table shown in Table~\ref{tab:ktab_lapqdpgradgrad}.
Note that the kernels from S points are not written in the combined kernel form \(L\oplus L^D\), because the charge S points and dipole S points do not have to be the same set of points.
Instead, for each source we use \(L\) if it is a charge and \(L^D\) if it is a dipole.
Because only charges are placed on M and L points, no $L^D$ kernel is necessary on the \MLtoML{} kernels in Table~\ref{tab:ktab_lapqdpgradgrad}.
This is possible because placing charges on M and L points is sufficient to match the field generated by dipole sources accurately.
This is how PVFMM implements this kernel summation problem \cite{MalhotraBiros2015PVFMMParallelKernel}.
This idea also works in many similar cases, such as \([q,\bQ]\to\left[\phi,\nabla\phi,\nabla\nabla\phi\right]\), or \([q,\bd,\bQ]\to\left[\phi,\nabla\phi,\nabla\nabla\phi\right]\).
\begin{table}[htbp]
  \centering
  \begin{tabular}{c|c|c|c}
    \hline
      & M                 & L                 & T                                                                                             \\
    \hline
    S & \(L\) and \(L^D\) & \(L\) and \(L^D\) & \(L\oplus \nabla L\oplus \nabla\nabla L\) and \(L^D\oplus \nabla L^D\oplus \nabla\nabla L^D\) \\
    M & \(L\)             & \(L\)             & \(L\oplus \nabla L\oplus \nabla\nabla L\)                                                     \\
    L & --                & \(L\)             & \(L\oplus \nabla L\oplus \nabla\nabla L\)                                                     \\
    \hline
  \end{tabular}
  \caption{\label{tab:ktab_lapqdpgradgrad} The kernel table to evaluate \([q,\bd]\to\left[\phi,\nabla\phi,\nabla\nabla\phi\right]\). }
\end{table}

\subsection{\label{subsec:stksl} Stokeslet and regularized Stokeslet}
The Stokeslet,
the primary Green's function for Stokes flow,
is a linear kernel \(\bG\) with dimensions \(3\times 3\):
\begin{align}
  G_{ij}=\frac{1}{8\pi}\left[\frac{\delta_{ij}}{r} +\frac{r_ir_j}{r^3} \right],
\end{align}
where \(\br=\bx-\by\), and we set the viscosity \(\mu=1\) for simplicity.
This kernel computes the velocity field \(\bu\) generated by a point force \(\bff\): \(u_i=G_{ij}f_j\).
In boundary integral methods, the Stokeslet is also called the Stokes single-layer kernel.
The kernel summation problem \(\bff\to\bu\) can be evaluated by using \(\bG\) directly in KIFMM.

In many applications, however, the Rotne-Prager-Yamakawa (RPY) kernel \cite{RotnePrager1969VariationalTreatmentHydrodynamic,Yamakawa1970TransportPropertiesPolymer,WajnrybMizerskiEtAl2013GeneralizationRotnePragerYamakawaMobility,ZukWajnrybEtAl2014RotnePragerYamakawa,MizerskiWajnrybEtAl2014RotnePragerYamakawa} is preferred because it approximates the hydrodynamic interactions among spheres with finite radii and guarantees that the mapping from \(\bff\) to \(\bu\) for all spheres is symmetric-positive-definite.
In other words, the energy is guaranteed to be dissipated by viscosity of the fluid, which fails in simulations that model the hydrodynamic interaction by the Stokeslet when the spheres get close to each other.

The RPY theory first computes the flow field \(\bu^\alpha \) at an arbitrary location \(\bx^\alpha \) generated by particles at S points \(\by^\beta \) as
\begin{align}\label{eq:RPYflow}
  \bu^\alpha  =\sum_{\beta\leq N_S} \bG^{RPY}\left(\bx^\alpha,\by^\beta, b^\beta\right)\bff^\beta,
\end{align}
where \(\bff^\beta \) is the force applied to the source sphere \(\beta \) with radius \(b^\beta \).
The RPY kernel \(\bG^{RPY}\) is a nonlinear function of \(b\),
\begin{align}\label{eq:RPYkernel}
  \bG^{RPY} = \left(1+\frac{1}{6}b^2\nabla^2\right) \bG,
\end{align}
where \(\bG\) is the Stokeslet.
The velocity of each target sphere \(\alpha\) at location \(\bx^\alpha\) is then constructed from the flow field \(\bu\) and its Laplacian \(\nabla^2 \bu\), to account for the sphere's finite radius \(a^\alpha\):
\begin{align}\label{eq:RPYtrgvel}
  \bu_a^\alpha = \left(1+\frac{1}{6}(a^{\alpha} )^2\nabla^2\right)\bu^\alpha.
\end{align}
If the radius \(a=0\), the target sphere is simply a tracer and its velocity equals that of the flow field.
Thus the RPY theory poses the kernel summation problem
\( [\bff,b]\to\left[\bu,\nabla^2\bu\right]. \)
After \(\bu\) and \(\nabla^2\bu\) are evaluated, the velocity of each target particle can be assembled using Eq.~(\ref{eq:RPYtrgvel}).
In other words, at the kernel summation stage, no information about the target particle radius \(a\) is needed.

The tricky part of this problem is the nonlinear dependence of \(\bG^{RPY}\) on \( b \).
If \(b\) is the same for all particles, \(b\) can be treated as a common parameter and \(\bG^{RPY}\) can still be written as a linear operator applied to \(\bff\).
In this particular case, the usual KIFMM algorithm still works.
However, in polydisperse cases the radius \(b\) is different for each particle and \(\bG^{RPY}\) cannot be used as \MLtoML{} kernels.

This can be conveniently accommodated with KAFMM.
We can continue to use the linear singular Stokeslet \(\bG\) for \MLtoML{} kernels because the velocity field generated by \(\bG^{RPY}\) as given by Eq.~(\ref{eq:RPYflow}) is still an analytical solution of the Stokes equation and can be matched by \(\bG\) with point-force sources without the finite radius $b$ on the M and L points.
Thus we can use the kernel table shown in Table~\ref{tab:ktab_rpy} for the RPY problem.
Summarizing how this table works: Source points S require the nonlinear \(\bG^{RPY}\) because of their finite radii.
Target points T require \(\bG\oplus \nabla^2 \bG\) to compute the effect of their finite radii.
The linear Stokeslet \(\bG\) suffices for the \MLtoML{} kernels.
\begin{table}[htbp]
  \centering
  \begin{tabular}{c|c|c|c}
    \hline
      & M             & L             & T                                      \\
    \hline
    S & \(\bG^{RPY}\) & \(\bG^{RPY}\) & \(\bG^{RPY}\oplus \nabla^2 \bG^{RPY}\) \\
    M & \(\bG\)       & \(\bG\)       & \(\bG\oplus \nabla^2 \bG\)             \\
    L & --            & \(\bG\)       & \(\bG\oplus \nabla^2 \bG\)             \\
    \hline
  \end{tabular}
  \caption{\label{tab:ktab_rpy} The kernel table to evaluate the RPY flow field \(\bu\) and its Laplacian: \(\left[\bff,b\right]\to\left[\bu,\nabla^2\bu\right]\).
  }
\end{table}

This is, of course, not the only way to compute the RPY kernel summation problem.
One straightforward approach using KIFMM is to split the problem into two sub-problems, utilizing the fact that $\bG^{RPY}\bff = \bG\bff + \nabla^2\bG \left(\frac{1}{6}b^2\bff\right)$ and the vanishing double Laplacian of Stokeslet $\nabla^2\nabla^2\bG=0$.
First sub-problem applies the combination kernel $\bG\oplus\nabla^2\bG$ on source point forces $\bff$ to compute part of $\bu$ and $\nabla^2\bu$, and the second sub-problem applies the kernel $\nabla^2\bG$ on the source point forces $\frac{1}{6}b^2\bff$ to compute the rest part of $\bu$.

Our method here is more compact, readily periodized reusing the M2L operator \cite{YanShelley2018Flexiblyimposingperiodicity} for $\bG$, and applicable to other cases where such splitting is not possible, for example, the regularized Stokeslet discussed in the following.

Above we have discussed only the case when the source and target particles do not overlap, i.e., \(\ABS{\br^{\alpha\beta}}>a^{\alpha}+b^{\beta}\).
If they overlap, the velocity of the target particle no longer satisfies Eq.~(\ref{eq:RPYtrgvel}) \cite{RotnePrager1969VariationalTreatmentHydrodynamic,ZukWajnrybEtAl2014RotnePragerYamakawa}, and it should be computed as an extra post-processing stage after FMM completes.
However, this post-processing stage involves corrections only for close pairs,
and with well-established near-neighbor detection algorithms it takes negligible time in comparison to the octree traversal.

The RPY kernel can also be thought of as a regularization of the singular Stokeslet \(\bG\), where the singular point force is regularized to a sphere of finite radius.
Another way to regularize \(\bG\) is to assume that the force is applied not at a single point but with a spatial distribution specified by a spreading function (or ``blob function'') \(\psi(\by-\by^\beta)\) around the source point \(\by^\beta\).
The spreading function may have infinite support but in general decays to zero at infinity and should satisfy \(\int \psi(\by)d^3\by=1\).
The singular Stokeslet can be considered the limiting case where \(\psi(\by-\by^\beta)=\delta(\by-\by^\beta)\).

One popular choice for \(\psi\) is: \cite{CortezFauciEtAl2005MethodRegularizedStokeslets,Cortez2018RegularizedStokesletSegments,OlsonLimEtAl2013ModelingDynamicsElastic}
\begin{align}
  \psi^\epsilon(\by)= \frac{15\epsilon^4}{8\pi\left(r^2+\epsilon^2\right)^{7/2}},
\end{align}
where \(\epsilon\) is the length scale of the spreading function \(\psi^\epsilon\).
\(\int \psi^\epsilon(\by)d^3\by = 1\) for any \(\epsilon >0\).
The flow field generated by this regularized force \(\psi^\epsilon \bff\) is \(\bu=\bG^\epsilon \bff\), where
\begin{align}
  G_{ij}^\epsilon = \frac{1}{8\pi}\left[\frac{r^{2}+2 \epsilon^{2}}{\left(r^{2}+\epsilon^{2}\right)^{3 / 2}} \delta_{i j} +\frac{1}{\left(r^{2}+\epsilon^{2}\right)^{3 / 2}} r_i r_j  \right].
\end{align}
The kernel table to compute \([\bff,\epsilon] \to \bu\) is similar to the RPY case without the Laplacian, as shown in Table~\ref{tab:ktab_reg}.
Similar to the case for \(\bG^{RPY}\),
our method allows different regularization parameters \(\epsilon\) to be used for different source points,
even though \(\bG^{\epsilon}\) depends on \(\epsilon\) nonlinearly.
\begin{table}[htbp]
  \centering
  \begin{tabular}{c|c|c|c}
    \hline
      & M                  & L                  & T                  \\
    \hline
    S & \(\bG^{\epsilon}\) & \(\bG^{\epsilon}\) & \(\bG^{\epsilon}\) \\
    M & \(\bG\)            & \(\bG\)            & \(\bG\)            \\
    L & --                 & \(\bG\)            & \(\bG\)            \\
    \hline
  \end{tabular}
  \caption{\label{tab:ktab_reg} The kernel table to evaluate the flow field \(\bu\) generated by a regularized force \(\bff\psi^\epsilon\): \(\left[\bff,\epsilon\right]\to\bu\). }
\end{table}

Our method for \(\bG^{\epsilon}\) can be further extended to include torque \(\bt\) at source points and angular velocity \(\bomega\) at target points.
That is, we evaluate the kernel summation problem \(\left[\bff,\bt,\epsilon\right]\to\left[\bu,\bomega\right]\).
We first define the necessary kernels in addition to \(\bG^{\epsilon}\)  \cite{OlsonLimEtAl2013ModelingDynamicsElastic}:
\begin{align}
  u_i      & =G_{ij}^\epsilon f_j,                                     & G_{ij}^{\epsilon}      & = H_1(r)\delta_{ij} + H_2(r)r_ir_j,                     \\
  \omega_i & =\Omega_{ij}^\epsilon f_j, \quad u_i=T_{ij}^\epsilon t_j, & \Omega_{ij}^{\epsilon} & =T_{ij}^{\epsilon}=\frac{1}{2}Q(r)\varepsilon_{ijk}r_k, \\
  \omega_i & =W_{ij}^\epsilon t_j,                                     & W_{ij}^{\epsilon}      & =\frac{1}{4}D_1(r)\delta_{ij} +\frac{1}{4}D_2(r)r_ir_j,
\end{align}
where \(\varepsilon_{ijk}\) is the Levi-Civita tensor and
\begin{align}
  H_{1}(r) & =\frac{2 \epsilon^{2}+r^{2}}{8 \pi\left(\epsilon^{2}+r^{2}\right)^{3 / 2}},                         \\
  H_{2}(r) & =\frac{1}{8 \pi\left(\epsilon^{2}+r^{2}\right)^{3 / 2}},                                            \\
  Q(r)     & =\frac{5 \epsilon^{2}+2 r^{2}}{8 \pi\left(\epsilon^{2}+r^{2}\right)^{5 / 2}},                       \\
  D_{1}(r) & =\frac{10 \epsilon^{4}-7 \epsilon^{2} r^{2}-2 r^{4}}{8 \pi\left(\epsilon^{2}+r^{2}\right)^{7 / 2}}, \\
  D_{2}(r) & =\frac{21 \epsilon^{2}+6 r^{2}}{8 \pi\left(\epsilon^{2}+r^{2}\right)^{7 / 2}}.
\end{align}
In the limit \(\epsilon\to 0\), the kernel functions \(\bG^\epsilon, \bOmega^\epsilon, \bT^\epsilon, \bW^\epsilon\) degenerate to the singular point-force and -torque kernels \(\bG, \bOmega, \bT, \bW\).

Our method has no difficulty handling the case where the regularization parameter \(\epsilon\) is different for force and torque for each point, but in most applications each S point has only one spreading parameter \(\epsilon\) shared by both force and torque.
For that case,
the \ST\ kernel \(\bG^\epsilon\oplus\bOmega^\epsilon\oplus\bT^\epsilon\oplus\bW^\epsilon\) has dimensions \(7\times 6\),
which computes \(\left[\bff,\bt,\epsilon\right]\to\left[\bu,\bomega\right]\).
The resulting kernel table is only slightly more complicated than the previous case, as shown in Table~\ref{tab:ktab_regftvw}.

Note that the Stokeslet kernel \(\bG\) suffices for \MLtoML{} points, which match the flow field with singular point-force sources, without needing special torque sources or a computation of angular velocity.
The other entries in the table can be easily understood as follows:
Since S points have regularized forces and torques, all kernels from S must use the regularized kernel functions.
Since we need \(\left[\bu,\bomega\right]\) on T points, the \ST\ kernel must evaluate those using the full regularized kernel combination.
The \MT\ and \LT\ kernels are the singular \(\bG\oplus\bOmega\) combination kernel because there is neither regularization nor torque sources at the M and L equivalent points.
\begin{table}[htbp]
  \centering
  \begin{tabular}{c|c|c|c}
    \hline
      & M                                    & L                                      & T                                                                          \\
    \hline
    S & \(\bG^{\epsilon}\oplus\bT^\epsilon\) & \(\bG^{\epsilon}\oplus\bT^{\epsilon}\) & \(\bG^\epsilon\oplus\bOmega^\epsilon\oplus\bT^\epsilon\oplus\bW^\epsilon\) \\
    M & \(\bG\)                              & \(\bG\)                                & \(\bG\oplus\bOmega\)                                                       \\
    L & --                                   & \(\bG\)                                & \(\bG\oplus\bOmega\)                                                       \\
    \hline
  \end{tabular}
  \caption{\label{tab:ktab_regftvw} The kernel table to evaluate \(\left[\bff,\bt,\epsilon\right]\to\left[\bu,\bomega\right]\). }
\end{table}

For summation problems (such as this one) that use a spreading function \(\psi^\epsilon\), the convergence error analysis results of the KIFMM algorithm \cite{YingBirosEtAl2004kernelindependentadaptive} imposes a limit on the magnitude of the spreading parameter \(\epsilon\).
A detailed quantitative proof of this subtlety is beyond the scope of this paper, but we offer a qualitative explanation:
The equivalence of the S and M points corresponding to an octree leaf box requires that the M points completely surround the S points.
If the spreading parameter \(\epsilon \) is too large, the regularized force or torque may significantly spread outside the range formed by the M points.
In this case, the field generated by the M points is no longer guaranteed to match the field generated by S points at locations faraway.
In other words, some ``faraway'' points are no longer well separated from the S points and the field generated by the S points can no longer be approximated by a small number of M points.
Technically, this subtlety requires that the characteristic spreading width \(\epsilon\) must be far smaller than the smallest leaf box in the adaptive octree.
This is typically true for practical problems, because \(\epsilon\) is usually chosen to be smaller than or comparable to the minimal point-point distance of S points,
while the smallest octree box is set to contain \(O(10^3)\) S points in order to maintain high efficiency of octree traversal.
In Section~\ref{subsec:open} we present numerical tests where this requirement on \(\epsilon\) is satisfied and we show that the convergence is not affected.

\subsection{\label{subsec:stkdl} Extensions to Stokeslets: pressure, source, sink, and double layer}
The Stokeslet and its regularized versions are widely used, but in many applications they are not sufficient.
First, the fluid pressure is not computed by \(\bG\), \(\bG^{RPY}\), or \(\bG^\epsilon\).
Second, if a source or sink of fluid appears somewhere, the flow field cannot be described by Stokeslets or their regularized versions anymore, because the divergence \(\nabla\cdot\bu = \nabla\cdot\left(\bG\bff\right)\) always vanishes.
Third, in boundary integral methods there are so-called Stokes double-layer operators and Stokes traction operators that depend on higher order kernels instead of Stokeslets.
In this section, we first extend the Stokeslet \(\bG\) to a new kernel \(\bG^P\) with dimensions \(4\times 4\) to handle the pressure field and sources/sinks of fluid.
We call \(\bG^P\) the ``StokesPVel'' kernel because it computes the pressure and velocity \([p,\bu]\) from the force and fluid sources/sinks \([\bff,q]\).
Next, we show how we implement the Stokes doublet and stress kernels used in boundary integral Stokes double-layer and traction operators using this StokesPVel kernel \(\bG^P\).

The kernel \(\bG^P\) is a linear operator applied to $\left[\bff,q\right]$, extended from the Stokeslet \(\bG\),
\begin{align}
  \bG^P(\br) = \begin{bmatrix}
                 \dfrac{1}{4\pi r^3}\br^T & 0                        \\
                 \bG(\br)                 & -\dfrac{1}{8\pi r^3} \br
               \end{bmatrix},
\end{align}
where \(\br=\bx-\by\in\bbR^3\) is a column vector and its transpose \(\br^T\) is a row vector.
The linear mapping \(\left[p,\bu\right]=\bG^P\left[\bff,q\right]\) can be written as
\begin{align}
  p   & =\dfrac{1}{4\pi r^3} r_j f_j,
  \label{eq:stkslp}                                                                                               \\
  u_i & = \frac{1}{8\pi}\left(\frac{r^2\delta_{ij}+r_ir_j}{r^3} f_j -\frac{r_i}{r^3}q\right). \label{eq:stkslvel}
\end{align}
Here we have defined the term \(q\) to be twice the actual source pumped into the fluid (or sink extracted) to obtain the convenient common \({1}/{8\pi}\) prefactor.
Also, the source or sink term does not contribute to fluid pressure \(p\).

To evaluate the kernel summation problem \(\left[\bff,q\right]\to\left[p,\bu\right]\), KIFMM can be used directly with the StokesPVel kernel \(\bG^P \) for all 8 interactions.
In applications requiring \(\nabla p, \nabla u , \) or \(\nabla^2 u\) on target points, the kernel tables are shown in
Tables~\ref{tab:ktab_stkslpvelgrad}.

\begin{table}[htbp]
  \centering
  \begin{tabular}{c|c|c|c}
    \hline
      & M         & L         & T                                                          \\
    \hline
    S & \(\bG^P\) & \(\bG^P\) & \(\bG^P\oplus\nabla\bG^P\) or \(\bG^P\oplus\nabla^2\bG^P\) \\
    M & \(\bG^P\) & \(\bG^P\) & \(\bG^P\oplus\nabla\bG^P\) or \(\bG^P\oplus\nabla^2\bG^P\) \\
    L & --        & \(\bG^P\) & \(\bG^P\oplus\nabla\bG^P\) or \(\bG^P\oplus\nabla^2\bG^P\) \\
    \hline
  \end{tabular}
  \caption{\label{tab:ktab_stkslpvelgrad} The kernel table to evaluate \([\bff,q]\to\left[p,\bu, \nabla p, \nabla\bu\right]\) or \([\bff,q]\to\left[p, \bu, \nabla^2 \bu \right]\). }
\end{table}

In boundary integral methods, the Stokes double-layer kernel is also widely used.
It's defined as a third-order tensor,
\begin{align}
  \Tcal_{ijk}= \frac{3}{8\pi}\left[ - \frac{r_ir_jr_k}{r^5} \right],
\end{align}
which is used in both the Stokes double-layer operator and the Stokes traction operator.
In the double-layer operator, it computes the flow field generated by an arbitrary tensorial source term \(D_{jk}\):
\begin{align}
  p   & =\frac{1}{4\pi}\left(-3\frac{r_jr_k}{r^5} + \frac{\delta_{jk}}{r^3}\right)D_{jk}, \label{eq:stkdlp} \\
  u_i & = \Tcal_{ijk} D_{jk}. \label{eq:stkdlvel}
\end{align}
For simplicity we write this as a linear mapping
\begin{align}
  \begin{bmatrix}
    p \\
    \bu
  \end{bmatrix} = \bK \bD.
\end{align}
In our implementation, the tensor \(\bD\) is flattened as a 9-dimensional vector and \(\bK\) has kernel dimensions \(9\times 4\).
This is more general than the case in boundary integral methods, where \(D_{jk}\) has only 6 degrees of freedom because it is constructed from a vectorial field \(v_j\) and a surface normal vector \(n_k\) as \(D_{jk}=v_jn_k\).

The flow field given by Eq.~(\ref{eq:stkdlvel}) actually contains a scalar source/sink term.
To see this, consider a sphere \(\Omega\) with radius \(R\) and place a tensorial source \(D_{jk}\) at the center.
The net flow \(q_{\scriptscriptstyle\Omega}\) out of the spherical surface is
\begin{align}
  q_{_\Omega}=\int_\Omega \bu\cdot\bn dS = -\frac{1}{2}\delta_{jk}D_{jk},
\end{align}
which is independent of the radius R.
The net flow \(q_{\scriptscriptstyle\Omega}\) is nonzero unless the source \(D_{jk}\) is trace-free.
This net flow cannot be matched by Stokeslets alone, so we use the StokesPVel kernel \(\bG^P\) and include scalar fluid sources \(q\) on M and L points during the octree traversal.
Full tensor sources \(D_{jk}\) are not needed, analogous to the Laplace summation problem with dipole or quadrupole sources on S points but only scalar charges on M and L points.

To evaluate the kernel summation problem \(\bD\to\left[p,\bu, \nabla p, \nabla \bu\right]\) or \(\bD\to\left[p, \bu, \nabla^2 \bu\right]\) we use the kernel table shown in Table~\ref{tab:ktab_stkdlpvelgrad},
where the combined kernel \(\bK\oplus\nabla\bK\) has kernel dimensions \(9\times 16\) and the combined kernel \(\bK\oplus\nabla^2\bK\) has kernel dimensions \(9\times 7\).
(The kernel \(\nabla^2\bK\) has kernel dimensions of only \(9\times 3\) because \(\nabla^2p\) is identically zero everywhere in Stokes flow and need not be computed.)

\begin{table}[htbp]
  \centering
  \begin{tabular}{c|c|c|c}
    \hline
      & M          & L         & T                                                           \\
    \hline
    S & \(\bK \)   & \(\bK \)  & \(\bK \oplus \nabla \bK \)  or  \(\bK \oplus \nabla^2 \bK\) \\
    M & \(\bG^P \) & \(\bG^P\) & \(\bG^P\oplus\nabla\bG^P \) or \(\bG^P\oplus\nabla^2\bG^P\) \\
    L & --         & \(\bG^P\) & \(\bG^P\oplus\nabla\bG^P \) or \(\bG^P\oplus\nabla^2\bG^P\) \\
    \hline
  \end{tabular}
  \caption{\label{tab:ktab_stkdlpvelgrad} The kernel table to evaluate \(\bD\to\left[p,\bu, \nabla p, \nabla \bu\right]\) or \(\bD\to\left[p, \bu, \nabla^2 \bu\right]\). }
\end{table}

Besides the double-layer operator, boundary integral methods also use the traction operator \(\bsigma\), which computes the fluid stress:
\begin{align}\label{eq:stressdef}
  \bsigma = -p\bI + \left[\nabla \bu + \left(\nabla\bu\right)^T\right],
\end{align}
where \(\bI\) is the identity tensor.
The fluid viscosity is set to 1.
In our implementation we write \(\sigma_{ij}\) as a flattened 9-dimensional vector.
Following our approach with the StokesPVel kernel, we include both a point force \(\bff\) and a source/sink term \(q\) on source points,
and we write the linear mapping \(\left[\bff,q\right]\to\bsigma\) as
\begin{align}
  \bsigma = \bS \begin{bmatrix}
                  \bff \\
                  q
                \end{bmatrix}.
\end{align}
This kernel function \(\bS\) has kernel dimensions \(4\times 9\).
Its entries can be determined by substituting the flow field of Eqs.~(\ref{eq:stkslp}) and (\ref{eq:stkslvel}) into the definition Eq.~(\ref{eq:stressdef}).

In the special case where the source/sink term \(q\) vanishes, the mapping \(\bff\to\bsigma\) invokes only the kernel \(\Tcal_{ijk}\):
\begin{align}
  \sigma_{ij} = \frac{3}{4\pi}\left(-\frac{r_ir_jr_k}{r^5}\right)f_k = 2 \Tcal_{ijk} f_k.
\end{align}
To evaluate the general kernel summation problem \([\bff, q ] \to \bsigma\), we construct the kernel table shown in Table~\ref{tab:ktab_stkstress}.
The kernel \(\bS\) is invoked only for target-point interactions because we only need the stress on target points.
For all other interactions, we only need to match the pressure and velocity fields on check surfaces, for which \(\bG^P \) suffices.
The definition Eq.~(\ref{eq:stressdef}) guarantees that the stress field must be equivalent if \(p\) and \(\bu\) are both equivalent.
\begin{table}[htbp]
  \centering
  \begin{tabular}{c|c|c|c}
    \hline
      & M         & L          & T        \\
    \hline
    S & \(\bG^P\) & \(\bG^P\)  & \(\bS\)  \\
    M & \(\bG^P\) & \(\bG^P\)  & \(\bS\)  \\
    L & --        & \(\bG^P \) & \(\bS \) \\
    \hline
  \end{tabular}
  \caption{\label{tab:ktab_stkstress} The kernel table to evaluate \([\bff,q]\to \bsigma\). }
\end{table}

Further, we consider the stress generated by the double-layer source \(D_{jk}\).
We write the mapping as \(\bsigma = \bS^D \bD\).
The kernel \(\bS^D\) has kernel dimensions \(9\times 9\), since we implement both \(\bsigma\) and \(\bD\) as flattened 9-dimensional vectors.
The entries of \(\bS^D\) can be conveniently computed with the flow field of Eqs.~(\ref{eq:stkdlp}) and (\ref{eq:stkdlvel}).
Similar to the previous cases, we place only \(\bff,q\) sources on M and L points, and use \(\bG^P\) for \MLtoML{} interactions.
The kernel table to evaluate \(\bD \to \bsigma\) is shown in Table~\ref{tab:ktab_stkdlstress}.
\begin{table}[htbp]
  \centering
  \begin{tabular}{c|c|c|c}
    \hline
      & M         & L          & T         \\
    \hline
    S & \(\bK\)   & \(\bK\)    & \(\bS^D\) \\
    M & \(\bG^P\) & \(\bG^P\)  & \(\bS\)   \\
    L & --        & \(\bG^P \) & \(\bS \)  \\
    \hline
  \end{tabular}
  \caption{\label{tab:ktab_stkdlstress} The kernel table to evaluate \(\bD\to \bsigma\). }
\end{table}

Last but not least, if both single-layer sources \(\bff,q\) and double-layer sources \(\bD\) exist in the system, we can evaluate the kernel summation problem by traversing the KAFMM octree once, in the same way as shown in Table~\ref{tab:ktab_lapqdpgradgrad}.
This is enabled by the flexible design of the software package PVFMM \cite{MalhotraBiros2015PVFMMParallelKernel}.

In this section we have frequently used the terms ``single layer'' and ``double layer'', borrowed from boundary integral methods.
This doesn't mean the KAFMM can only be applied to boundary integral methods.
For example, in Stokesian Dynamics or other singularity methods based on multipole expansion, the flow field (and its gradients and Laplacian) generated by a trace-free and symmetric stresslet \(S_{jk}\) is often computed, as given by
\begin{align}
  p   & =\frac{1}{4\pi}\left(-3\frac{r_jr_k}{r^5}\right)S_{jk},  \label{eq:stressletp}         \\
  u_i & = \frac{1}{8\pi}\left[ - \frac{3r_ir_jr_k}{r^5} \right]S_{jk}. \label{eq:stressletvel}
\end{align}
For this problem, our implementation for the arbitrary double-layer flow field of Eqs.~(\ref{eq:stkdlp}) and (\ref{eq:stkdlvel}) can be readily used without any modifications.

\section{\label{sec:bench}Implementation and Benchmarks}
As mentioned in the previous section, any kernel summation problems sharing the same \MLtoML{} kernels in their kernel tables can be evaluated together during a single pass of octree traversal.
For example, the potential generated by Laplace monopoles, dipoles, and quadrupoles can be computed with one tree traversal.
However, allowing combinations of an arbitrary number of such summation problems would significantly increase the complexity of the code base.
Our implementation of KAFMM is based on the high-performance KIFMM library \texttt{PVFMM} \cite{MalhotraBiros2015PVFMMParallelKernel}, which can handle a combination of two such problems.
Therefore, our package supports kernel summation problems for two sets of sources and one set of targets.
We call the two sets of sources ``single-layer'' (SL) sources and ``double-layer'' (DL) sources, and we always assign the higher order kernel function to DL sources.
For example, for a problem where both Laplace monopoles and dipoles exist, the DL sources are the dipoles.
Similarly, for a problem involving both StokesPVel sources \(\bff,q\) and Stokes double-layer sources \(\bD\), the DL sources refer to those \(\bD\) sources.
Our implementation allows the SL, DL, and target points to be three completely independent sets of points,
while making no assumptions about the numbers of points or their distributions.

We also implement singly, doubly, and triply periodic boundary conditions for all kernels using our flexible periodization method \cite{YanShelley2018Flexiblyimposingperiodicity}.
This periodization scheme is convenient because it relies on a precomputed translational operator \(\Tcal_{M2L}\), which adds the far-field effects from periodic images to the target points.
This precomputed operator depends on only the \ML\ kernel in the kernel table.
Therefore, the Laplace monopole, dipole, and quadrupole kernel functions all share one \(\Tcal_{M2L}\) operator;
the Stokeslet, regularized Stokeslet, and RPY kernel all share one operator;
and the StokesPVel kernel, Stokes double-layer and traction kernels and their derivatives all share one operator.
We have also implemented the no-slip wall boundary condition for Stokeslet and RPY kernels using our universal image method \cite{YanShelley2018UniversalImageSystems} with support for periodic boundary conditions.
See Table~\ref{tab:implement} in Appendix~\ref{app:stkfmm} for the full list of supported kernels and boundary conditions.
More software engineering details can be found in Appendix~\ref{app:simd}.

In this section we present numerical convergence and timing results for a few representative cases among the wide range of possible combinations of kernels and boundary conditions.
All these results were generated on one cluster node with two 20-core Intel Xeon Gold 6248 Skylake CPUs @ 2.5GHz.
We launched one MPI rank for each CPU socket with 20 OpenMP threads and bound each thread to one CPU core.
Hyper-threading was disabled, and the machines were set to ``performance'' mode to disable CPU idle down-clocking.
Detailed scaling tests is beyond the scope of this work because we did not modify any parallelization facilities in the \texttt{PVFMM} library.
Brief test results for large scale scaling benchmarks are incldued in Appendix~\ref{app:scaling}.

\begin{figure}[htbp]
  \includegraphics[width=\linewidth]{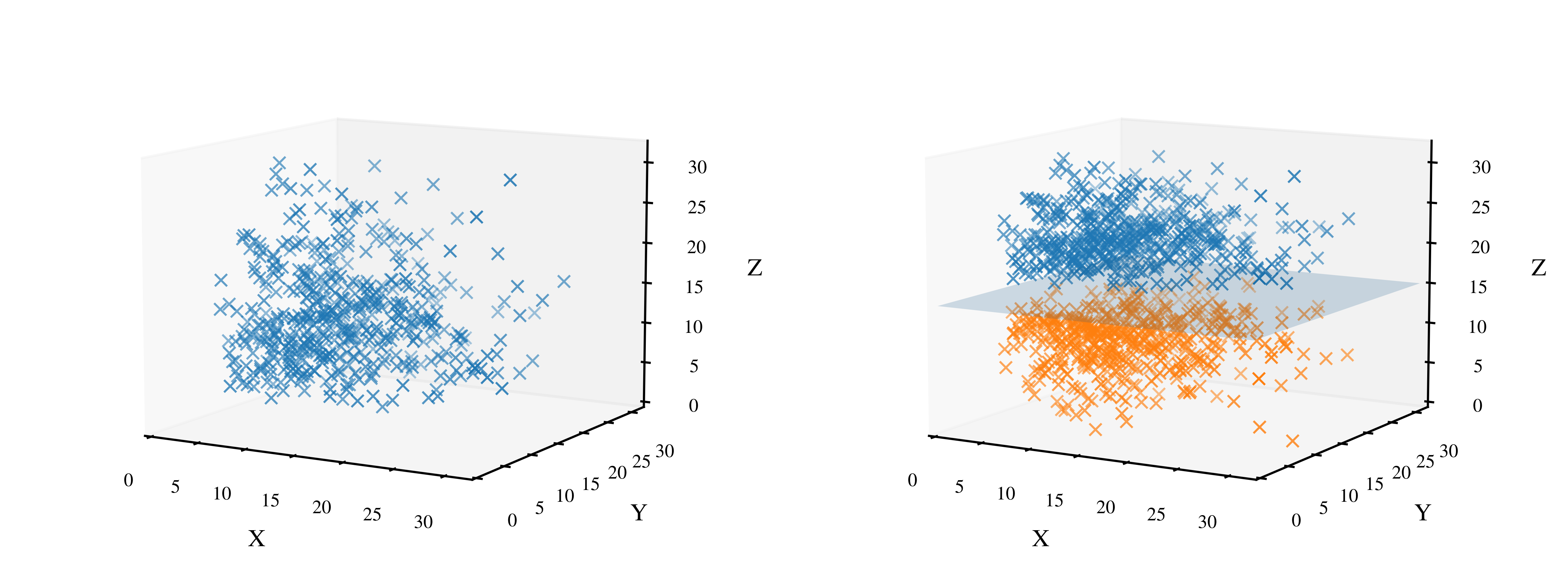}
  \caption{\label{fig:ptdist}
    The distribution of points in the cubic computational domain used in the convergence benchmarks for the cases without a no-slip wall (left) and with a no-slip wall at the \(z=15\) plane (right).
    The orange points below the no-slip wall are image points.
    Only 0.05\% of all points are shown here.
  }
\end{figure}

For simplicity, we fixed the SL, DL, and target points to be the same set of \(10^6\) points, as shown in Fig.~\ref{fig:ptdist},
with the points non-uniformly distributed in a cubic box defined by \([0,L)^3\).
The coordinates \((x,y,z)\) of each point were generated from the same random distribution,
\begin{align}\label{eq:ptdist}
  x\sim y\sim z \sim L (u-[u]),\quad p(u) = \frac{1}{s u \sqrt{2 \pi}} \exp \left(-\frac{(\ln u-m)^{2}}{2 s^{2}}\right),
\end{align}
where \(u\) is a random number following the log-normal distribution with PDF \(p(u)\), and \([u]\) is the maximum integer not larger than \(u\).
The parameters for the log-normal distribution were chosen to be \(m=-1,s=0.5\).
For the cases with a no-slip wall, the SL, DL, and T points and their images must both fit in the cubic box to construct the octree. We therefore put the wall at the \(z=L/2\) plane and shifted the z coordinate for each point to the upper half of the box:
\begin{align}
  z\sim \frac{L}{2} +\frac{L}{2}(u-[u]).
\end{align}
Their corresponding image sources were located in the lower half of the box, as shown in Fig.~\ref{fig:ptdist}.
We used \(L=32\) for all the benchmarks reported here.
Each component of the source values was generated from a random number with uniform distribution in the range \([-1,1]\), unless described otherwise.

The trade-off between computational time and accuracy in KIFMM is controlled by the multipole order \(m\).
The number of M and L points surrounding each leaf octree box is \(M=6(m-1)^2+2\).
The total computational cost of octree traversal is proportional to \(M^3\log M\) \cite{MalhotraBiros2015PVFMMParallelKernel}.
We set \(m_{\max}=16\), and evaluated the \(L_2\)-norm relative convergence error $\epsilon_{L2}$ of all target values relative to this \(m_{\max}=16\), starting from \(m=6\):
\begin{align}
  \epsilon_{L2} = \frac{\lVert{\bg_m-\bg_{m=16}}\rVert_2}{\lVert{\bg_{m=16}}\rVert_2}.
\end{align}

The adaptive octree was constructed by recursive refinement until each leaf octree box contained no more than \(2000\) target points.

\subsection{\label{subsec:open}Open boundary condition}
\subsubsection{Laplace kernels}
\begin{figure}[htbp]
  \centering
  \begin{subfigure}[t]{.48\textwidth}
    \centering
    \includegraphics[width=\linewidth]{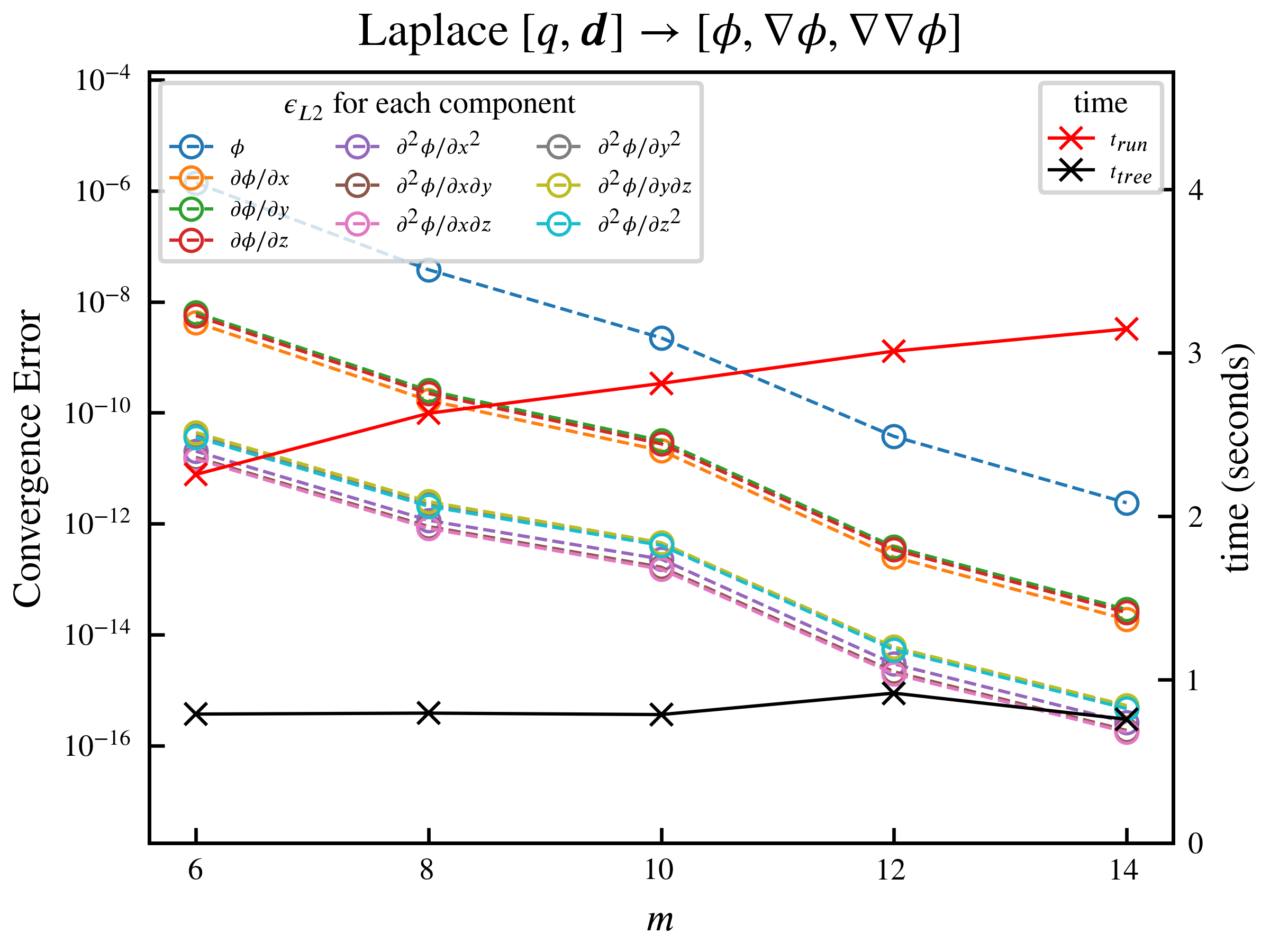}
    \caption{}
  \end{subfigure}
  \begin{subfigure}[t]{.48\textwidth}
    \centering
    \includegraphics[width=\linewidth]{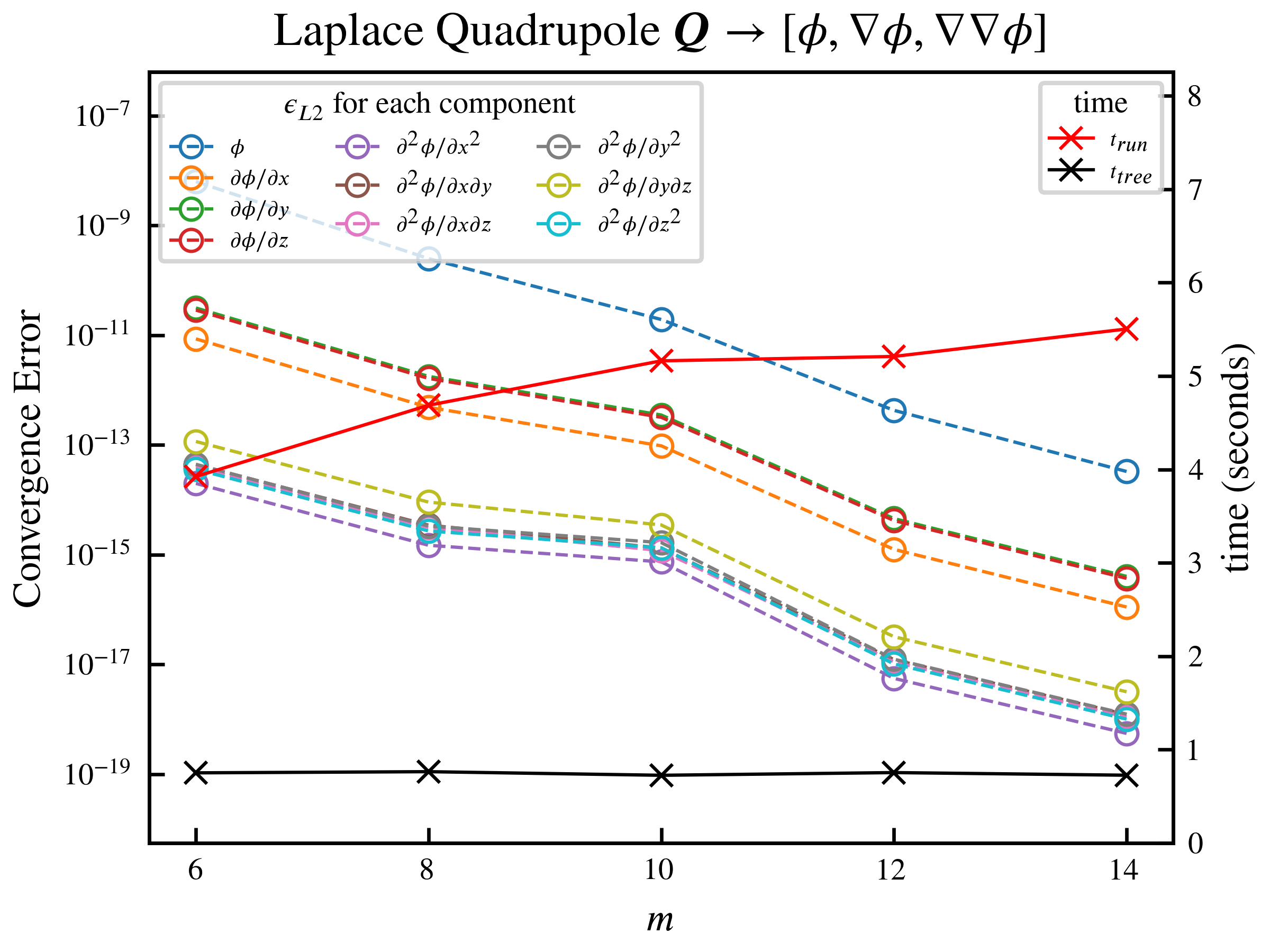}
    \caption{}
  \end{subfigure}
  \caption{\label{fig:lapP0N1}
    Convergence error and timing results for Laplace kernel problems with open boundary conditions for increasing multipole order \(m\) on octree boxes.
    (a) Results for the potential field \(\phi\) and its gradients generated by the combination of Laplace monopoles \(q\) on single-layer (SL) source points and dipoles \(\bd\) on double-layer (DL) source points (i.e., using kernel \texttt{LapPGradGrad} in the \texttt{STKFMM} package).
    (b) Results for \(\phi\) and its gradients generated by Laplace quadrupole sources \(\bQ\) (\texttt{LapQPGradGrad}).
    The convergence error is expressed as the \(L_2\)-norm error \(\epsilon_{\scriptscriptstyle L2}\) relative to target values computed using \(m=16\).
  }
\end{figure}

We benchmarked the performance with Laplace kernels and open boundary conditions by performing two kernel summation problems, which each computed \(\phi, \nabla\phi,\) and \(\nabla\nabla\phi\).
Fig.~\ref{fig:lapP0N1}(a) shows convergence and timing results for the first problem, in which the fields were generated by Laplace monopoles \(q\) on the SL points and dipoles \(\bd\) on the DL points.
Fig.~\ref{fig:lapP0N1}(b) shows the results for the second problem, which used Laplace quadrupoles \(\bQ\) on source points.
In both cases, as \(m\) linearly increases, the convergence error exponentially decreases for all components of the target values.
The time needed to construct the adaptive octree (\(t_\text{tree}\)) and then run the KAFMM algorithm by traversing the tree (\(t_\text{run}\)) both slightly increase with increasing multipole order \(m\).

\subsubsection{Stokes kernels}
\begin{figure}[htbp]
  \centering
  \begin{subfigure}{.48\textwidth}
    \centering
    \includegraphics[width=\linewidth]{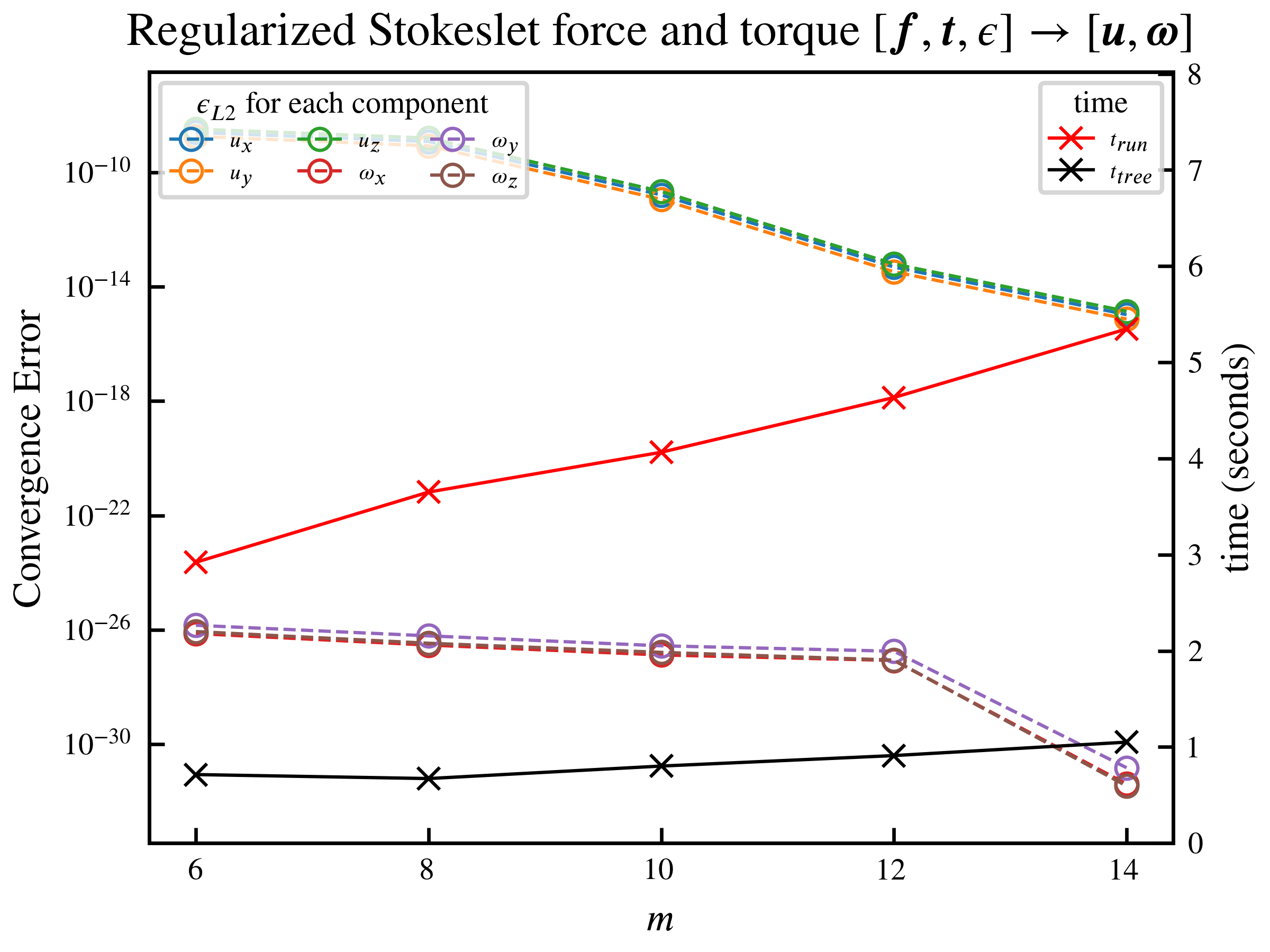}
    \caption{}
  \end{subfigure}
  \begin{subfigure}{.48\textwidth}
    \centering
    \includegraphics[width=\linewidth]{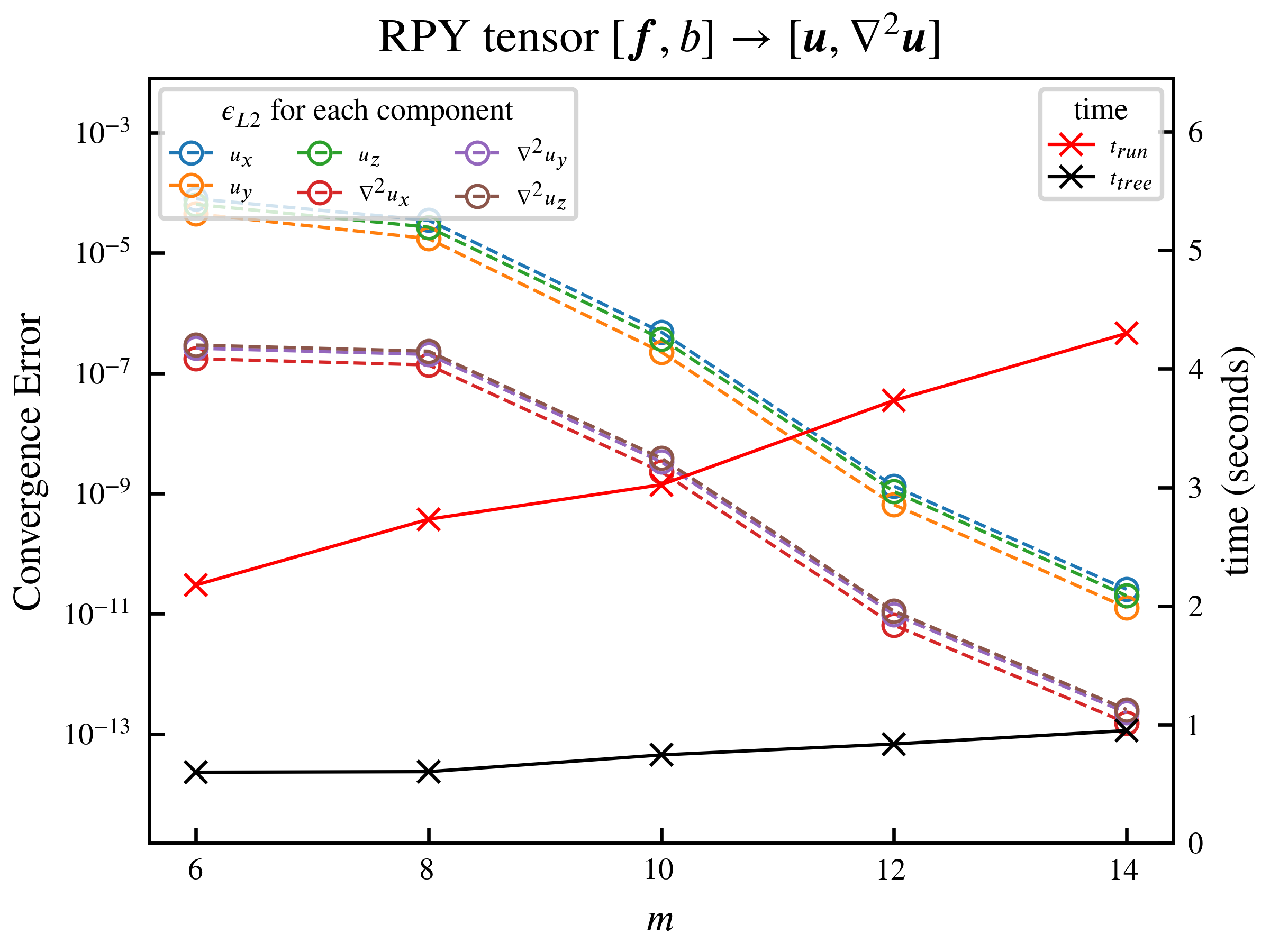}
    \caption{}
  \end{subfigure}
  \caption{\label{fig:stkregrpyP0N1}
    Convergence error and timing results for variants of the Stokeslet kernel.
    (a) Translational and rotational velocity fields \(\bu,\bomega\) generated by regularized Stokes force \(\bff\) and torque \(\bt\) sources (\texttt{StokesRegVelOmega}).
    (b) Velocity field and its Laplacian generated by RPY sources with radius \(b\)  (\texttt{RPY}). The regularization scale \(\epsilon\) and the radius \(b\)
    are randomly generated for each source point.
  }
\end{figure}

The Stokes kernels can be categorized into two groups.
The first group does not need to compute the pressure field \(p\), and uses the Stokeslet \(\bG\) as the \MLtoML{} kernels.
Two examples are the regularized Stokeslet and the RPY tensor.
Fig.~\ref{fig:stkregrpyP0N1} shows the benchmarks for these two kernels applied to summation problems in polydisperse systems.
Fig.~\ref{fig:stkregrpyP0N1}(a) shows the benchmarks for the regularized Stokeslet problem \([\bff,\bt,\epsilon]\to[\bu,\omega]\).
Fig.~\ref{fig:stkregrpyP0N1}(b) shows the benchmarks for the RPY problem \([\bff,b]\to[\bu,\nabla^2\bu]\).
The regularization scale \(\epsilon\) and the radius \(b\) were both randomly generated for each source point from a uniform distribution within the range \([0,10^{-4}]\).
Here \(10^{-4}\) was chosen to make sure that at the lowest level of the octree, the leaf box size was much larger than \(\epsilon\),
as is required by the convergence theorem of KIFMM.
For both kernels we again see exponentially decreasing convergence error and linearly increasing computational time as the order \(m\) increases.

\begin{figure}[htbp]
  \centering
  \begin{subfigure}{.48\textwidth}
    \centering
    \includegraphics[width=\linewidth]{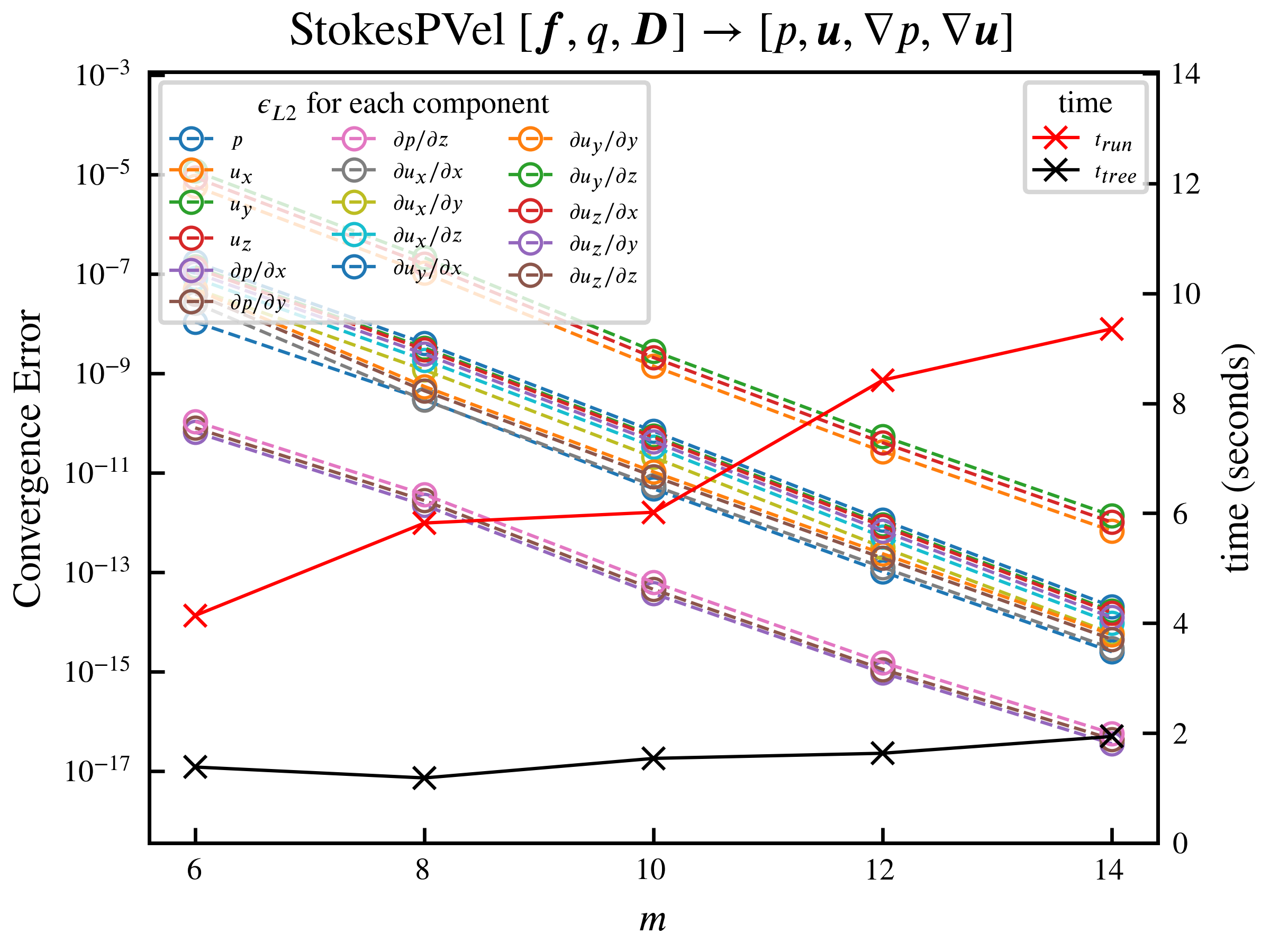}
    \caption{}
  \end{subfigure}
  \begin{subfigure}{.48\textwidth}
    \centering
    \includegraphics[width=\linewidth]{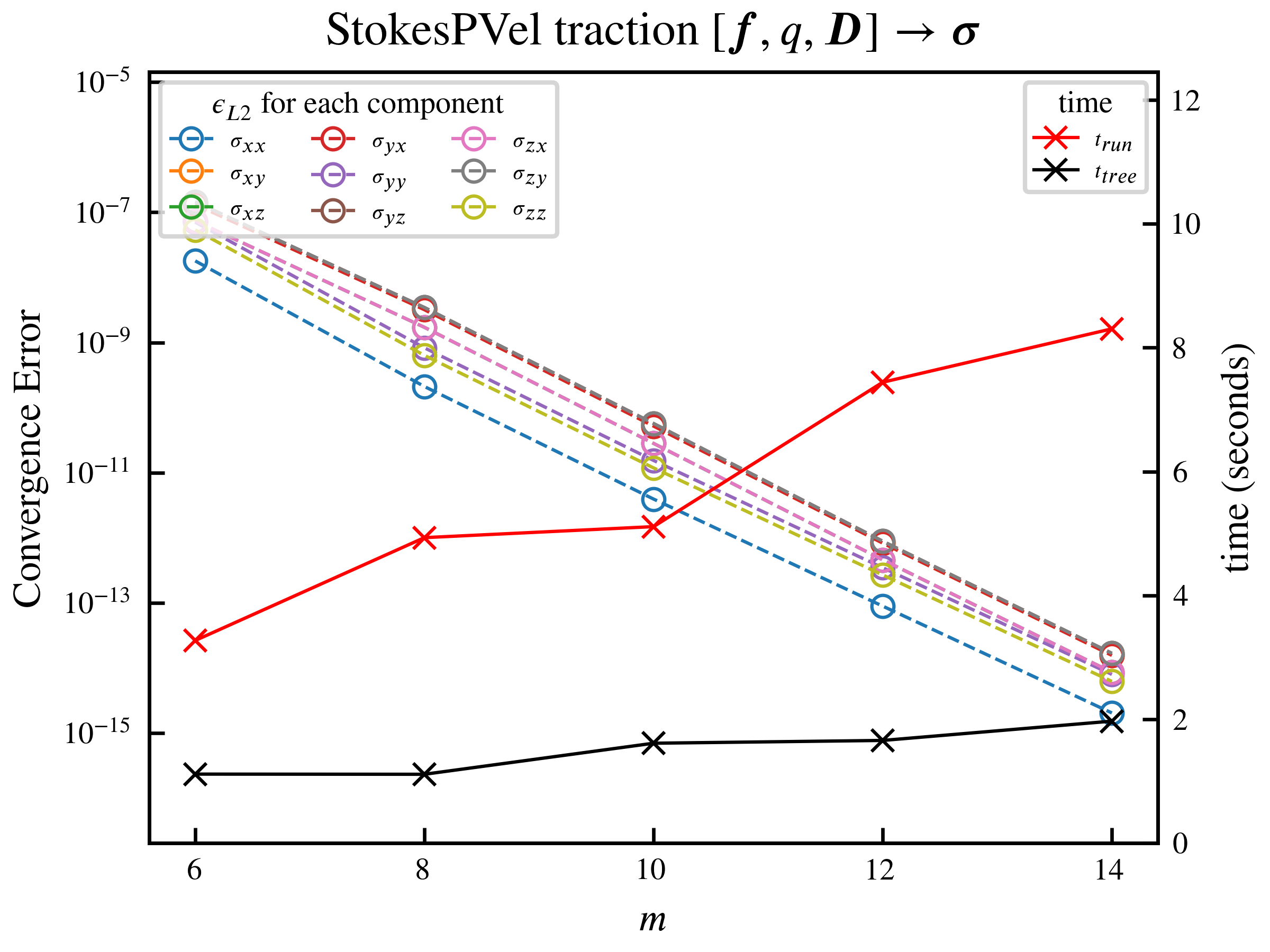}
    \caption{}
  \end{subfigure}
  \caption{\label{fig:stkpvelP0N1}
    Convergence error and timing results for combined StokesPVel and double-layer kernels.
    The StokesPVel kernel handles point forces \(\bff\) and fluid sources/sinks \(q\) on each SL point, while the double-layer kernel handles a double-layer source \(\bD\) on each DL point.
    Results for the computation of (a) the pressure \(p\) and velocity field \(\bu\) and their gradients (\texttt{PVelGrad}), and (b) the fluid stress \(\bsigma\)  (\texttt{Traction}).
  }
\end{figure}

The second group of Stokeslet kernels compute the pressure field and may have fluid sources or sinks in the system.
In this case we must use the extended StokesPVel kernel \(\bG^P\) as the \MLtoML{} kernels.
Fig.~\ref{fig:stkpvelP0N1} shows benchmarks for two kernel summation problems of this type.
In both cases the SL source points each have a point force \(\bff\) and a sink/source \(q\), and the DL points each have a double-layer source \(\bD\), flattened from a \(3\times 3\) tensor to a \(9\)-dimensional vector.
We show the results for computation of \([p,\bu,\nabla p,\nabla\bu]\) at each target point in  Fig.~\ref{fig:stkpvelP0N1}(a),
and for a computation of the hydrodynamic stress \(\bsigma\) at each target point in
Fig.~\ref{fig:stkpvelP0N1}(b).
Again, we see the same trend of exponentially improving accuracy and linearly increasing run times as the order \(m\) increases.

\subsection{\label{subsec:wall}No-slip Wall}
\begin{figure}[htbp]
  \centering
  \begin{subfigure}{.48\textwidth}
    \centering
    \includegraphics[width=\linewidth]{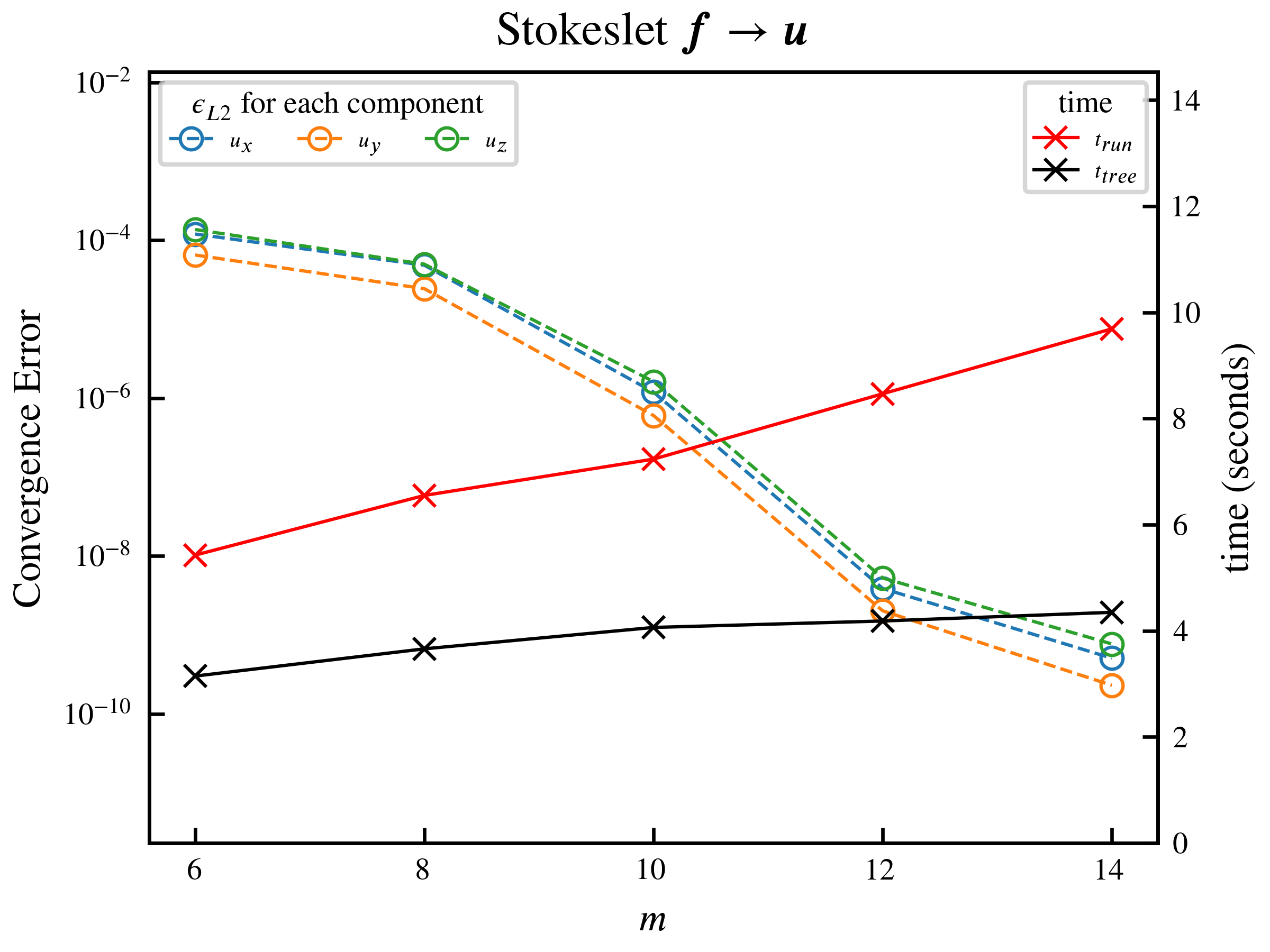}
    \caption{}
  \end{subfigure}
  \begin{subfigure}{.48\textwidth}
    \centering
    \includegraphics[width=\linewidth]{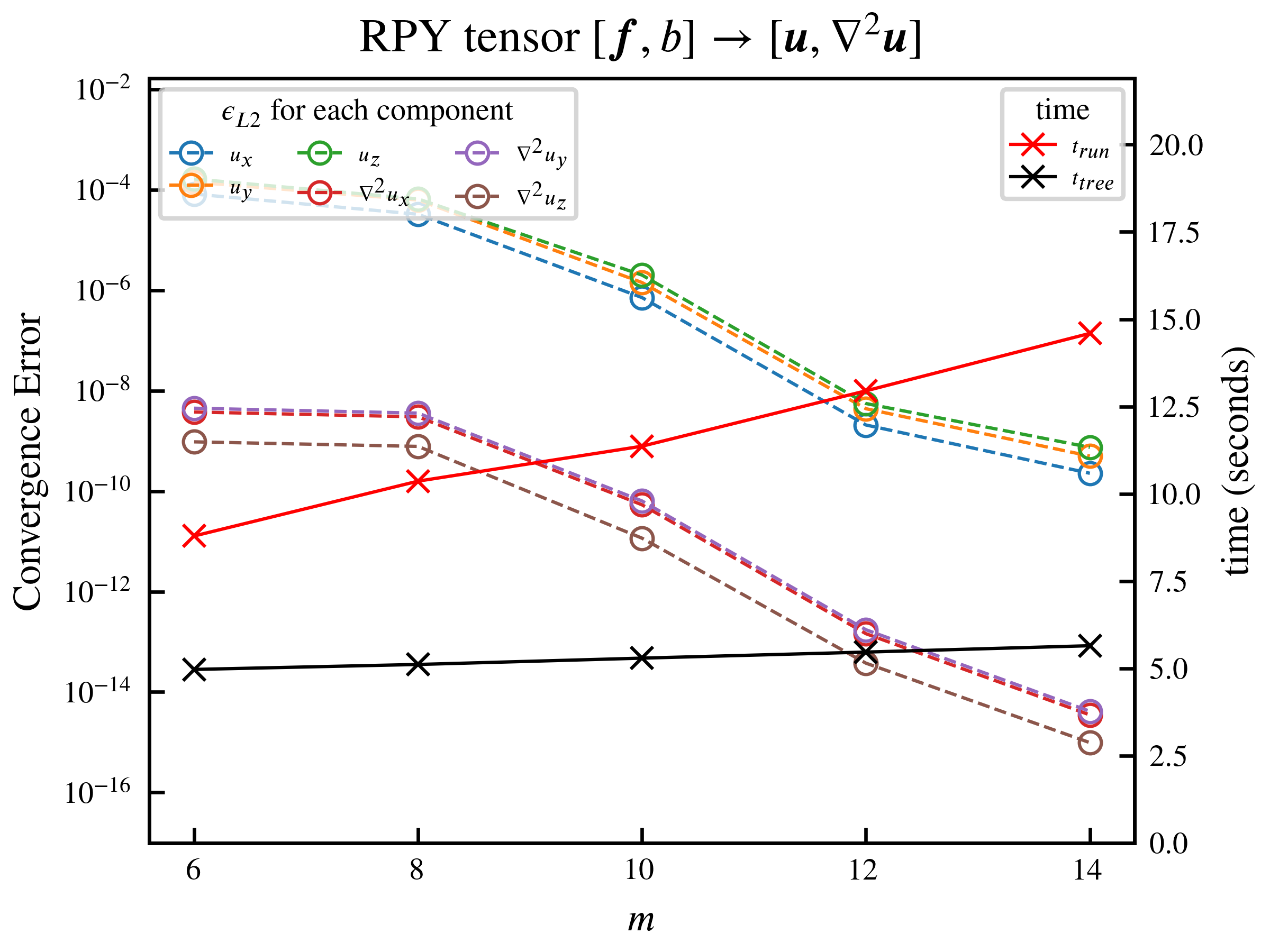}
    \caption{}
  \end{subfigure}
  \caption{
    \label{fig:wallP0N1}
    Convergence error and timing results for Stokeslet and RPY tensor summation problems with a no-slip wall and open boundary conditions.
    (a) Stokeslet computation of the velocity field generated by point force \(\bff\).
    (b) RPY tensor computation of the velocity field and its Laplacian generated by force \(\bff\) applied on particles with radius \(b\), where \(b\) is randomly generated for each particle.
  }
\end{figure}

Fig.~\ref{fig:wallP0N1} shows benchmarks for evaluating Stokeslets and RPY kernels above a no slip-wall, using the image method.
The image systems for such problems involve not only the Stokeslet or RPY sum with a set of image forces included, but also a couple of Laplace sums necessary to impose the no-slip condition, in addition to the image point forces.
These separate summations are computed by several independent calls of the KAFMM algorithm and then the results are summed together to generate the flow field above a no-slip wall.
Detailed derivation and summation can be found in our previous work \cite{YanShelley2018UniversalImageSystems}.

The summation of Stokeslets above a no-slip wall includes one Stokeslet summation and two Laplace monopole \& dipole summations.
The benchmark is shown in Fig.~\ref{fig:wallP0N1}(a).
The summation of RPY tensors above a no-slip wall is more complicated.
For a polydisperse system, 7 FMM evaluations were necessary in our previous work \cite{YanShelley2018UniversalImageSystems} where we used only Stokeslet and Laplace summations.
With KAFMM, we can efficiently implement the freespace RPY summation.
Therefore, we need only 4 summations, as tabulated in Table~\ref{tab:rpywall}.
Once $\bu^S,\phi^{SD},\phi^{SDZ},\phi^Q$ and their necessary derivaives have all been evaluated, the flow field at each target point $\bx$ can be assembled as
\begin{align}
  % trgValSDZ[4 * i + 1] + x3 * trgValSD[10 * i + 1] + x3 * trgValQ[10 * i + 1];
  % - trgValSD[10 * i] - trgValQ[10 * i]
  \bu         & =\bu^S+ \nabla \phi^{SDZ} + x_3 \nabla\phi^{SD} + x_3 \nabla\phi^Q - \left[0,0,\phi^{SD}+\phi^Q\right], \\
  \nabla^2\bu & =\nabla^2\bu^S + 2 \dpone{}{x_3}\nabla \left(\phi^{SD}+\phi^{Q}\right).
\end{align}
Here we used $1,2,3$ instead of $x,y,z$ to represent the coordinate components of target locations $\bx$ and source locations $\by$ to avoid confusion.
We assumed that the wall is on the plane $x_3=0$.

The benchmark for this case is shown in Fig.~\ref{fig:wallP0N1}(b).
Again, we see the same trend of exponentially decreasing errors and linearly increasing times when \(m\) increases.

\begin{table}[htbp]
  \centering
  \begin{tabular}{c|c|c|c}
    \hline
    source at $\by$                                   & source at $\by^I$                                  & summation       & target values                                     \\
    \hline
    $\bff_{12},b$                                     & $-\bff_{12},b$                                     & RPY             & $\bu^S,\nabla^2\bu^S$                             \\
    $f_3$                                             & $-f_3$ and $y_3(-f_1,-f_2,f_3)$                    & Laplace SL + DL & $\phi^{SD},\nabla\phi^{SD},\nabla\nabla\phi^{SD}$ \\
    $\frac{1}{2}y_3f_3$ and $\frac{1}{6}b^2(0,0,f_3)$ & $-\frac{1}{2}y_3f_3$ and $\frac{1}{6}b^2(0,0,f_3)$ & Laplace SL + DL & $\phi^{SDZ},\nabla\phi^{SDZ}$.                    \\
    $\frac{1}{3}b^2\begin{bmatrix}
                       f_3 & 0   & 0 \\
                       0   & f_3 & 0 \\
                       f_1 & f_2 & 0 \\
                     \end{bmatrix}$                    & 0                                                  & Laplace Q       & $\phi^Q,\nabla\phi^Q, \nabla\nabla\phi^Q$          \\
    \hline
  \end{tabular}
  \caption{\label{tab:rpywall}
    Four summations for the RPY tensor above no-slip wall problem.
    $\bx,\by$ represent the target and source locations.
    $\by^I$ refers to source image locations.
    All target values are evaluated at the target locations $\bx$.
    To avoid confusion, we use $1,2,3$ to represent the component of $\bx,\by$ and $\bff$.
    $\bff_{12}$ represent the components of the input force $\bff$ that is paralell with the no-slip wall.
    $b$ is the sphere radius of the source particles in the RPY problem.
  }
\end{table}

\subsection{\label{subsec:pbc}Periodic boundary conditions}
To benchmark the accuracy of periodic boundary condition, we compute the ``translation error'' in addition to the convergence error.
The translation error is defined as the \(L_2\)-norm errors of target value before and after all source and target points are translated by a random amount in the periodic directions.
For example, if the doubly periodic boundary condition is imposed on the \(x\)-\(y\) directions, we shift all source and target points by a randomly chosen vector \((r_x,r_y,0)\), modulo the periodicity
(i.e., any point that is shifted out of the cubic box is replaced by the periodic image that is shifted into the box).
If the periodic boundary condition is exactly satisfied, this translation will produce no errors in the target values.
For each \(m\) we pick one random translation and evaluate the target values before and after this translation to compute the translation error.

\begin{figure}[htbp]
  \centering
  \includegraphics[width=\linewidth]{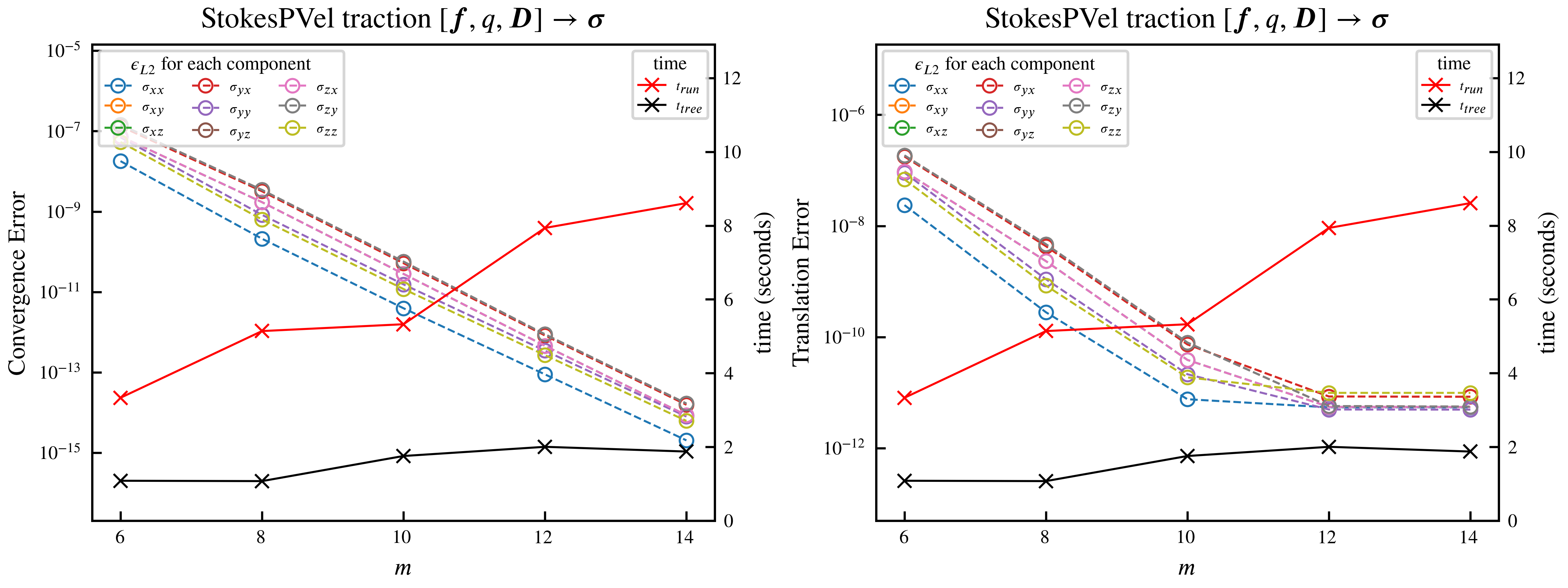}
  \caption{\label{fig:stktracP2N1}
    Convergence (left) and translation (right) errors and timing results for StokesPVel traction computations with doubly periodic boundary conditions.
    The fluid stress generated by StokesPVel single- and double-layer sources \(\bff,q, \bD \) was computed.
  }
\end{figure}

Fig.~\ref{fig:stktracP2N1} shows the timing results, convergence errors, and translation errors for the kernel summation problem \([\bff,q,\bD]\to\bsigma\) with doubly periodic boundary condition.
The convergence error shows the same exponentially decreasing trend when \(m\) increases as all previous cases.
However, the translation error may stagnate at \(10^{-11}\) in some cases.
We thoroughly tested this and found that this is induced by floating point round-off error.
More specifically, this happens when the minimal distance between two points get very close, e.g., in our non-uniformly distributed test points.
Our tests show that with random but uniformly distributed points, this problem is alleviated and may stagnate at \(10^{-12} - 10^{-13}\), while for a regular mesh of points this problem disappears.

This problem can be reproduced by simply evaluating the kernel function from only one source point to one target point, without any FMM operations.
Assume the two points are located at $\bx+\br$ and $\bx$.
After the translation along a vector $\bt$, which may be sufficiently large such that the two points are moved out of the periodic box, we have the error $\be=[(\bx+\br+\bt)-(\bx+\bt)] - [(\bx+\br)-\bx] \neq \bm{0}$ due to finite precision of floating point numbers.
When $\br$ is small, $\be$ becomes significant.
This error then enters the kernel evaluations because all kernels depends on the accuracy of $\br$.
These two-point tests show that this round-off relative L2-Norm error of kernel evaluation is around $10^{-12}$ when the distance $\br\approx 10^{-4}$, consistent with our FMM benchmark results for non-uniform points where the minimal distance between points is also about $10^{-4}$.

Unfortunately, there is no easy fix for this problem, except using extended precision floating point types (80 or even 128 bits).
However, in most applications there is usually a lower bound for the minimal distance between points, such as the grid point resolution or particle radius and this problem can be partially alleviated.

\begin{figure}[htbp]
  \includegraphics[width=\linewidth]{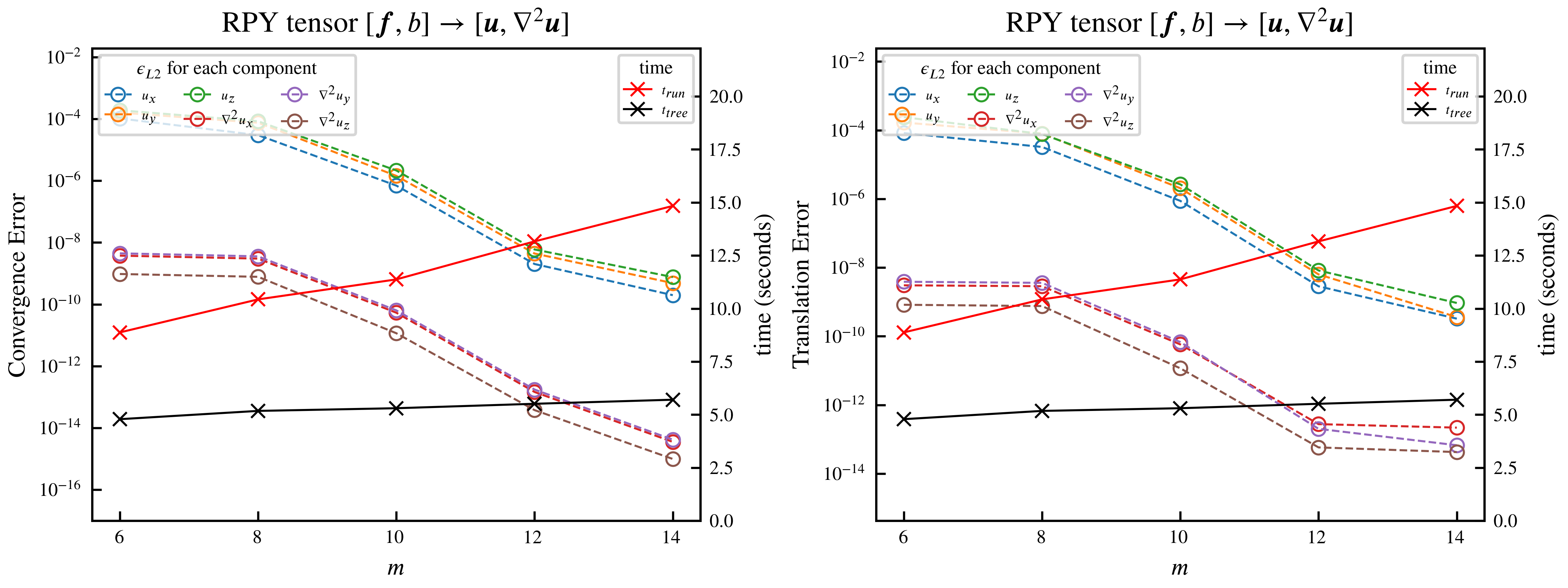}
  \caption{    \label{fig:wallrpyP1N1}
    Convergence (left) and translation (right) errors and timing results for a computation of the fluid velocity field and its Laplacian generated by RPY kernel sources \(\bff,b\) with singly periodic boundary condition, above a no-slip wall.
    The source values \(\bff,b\) are randomly generated for each source point.
  }
\end{figure}

Fig.~\ref{fig:wallrpyP1N1} shows the results for RPY tensor summation above a no-slip wall, with a singly periodic boundary condition imposed along the \(x\)-axis direction.
In this case the convergence errors and translation errors both converge to machine precision without issues.

The overhead of our method for implementing periodic boundary conditions is benchmarked in the Appendix~\ref{app:pbc}, where we compare timing for periodic and open boundary conditions for 4 different kernels.
Usually we have about 5\% - 10\% overhead, relative to the open boundary condition.

\section{\label{sec:conclusion}Conclusion}
In this work we extend the idea of `equivalent points' in classic KIFMM in a couple of kernels frequently used in the study of Stokes flow.
We have implemented this approach as an open-source package \texttt{STKFMM} to make it ready for use as a black-box---users do not need to know any of the implementation details.
A part of this implementation has been used in our previous simulation tasks \cite{LaGroneCortezEtAl2019ComplexDynamicsLong,YanCoronaEtAl2020ScalableComputationalPlatform} and here we presented a more systematic overview of the method.
The package \texttt{STKFMM} can also be easily extended to additional kernel summation problems with new kernel tables, thanks to the flexible design of the underlying PVFMM \cite{MalhotraBiros2015PVFMMParallelKernel} package.
For example, for the RPY tensor the package currently includes only the most popular force-velocity part,
but some applications require the mappings between force-angular velocity and torque-velocity\cite{WajnrybMizerskiEtAl2013GeneralizationRotnePragerYamakawaMobility,ZukWajnrybEtAl2014RotnePragerYamakawa}.
Those extensions can be conveniently incorporated, and periodic boundary conditions automatically work for those cases because they all rely on the same \MLtoML{} kernels as the $\bG^{RPY}$ kernel.

Since its invention more than a decade ago, the KIFMM method has been very successful thanks to its adaptivity and linear complexity.
We hope our open-source \texttt{STKFMM} package becomes a readily usable tool for researchers conducting simulations of Stokes flow in a wide range of areas such as rheology, biology, and chemical engineering.

\section*{Acknowledgement}
We thank Dhairya Malhotra and Alex Barnett for inspiring discussions, and thank Bryce Palmer for comments on the manuscript.

\appendix
\section{\label{app:stkfmm}All kernels supported by the open-source \texttt{STKFMM} package}

\begin{table}[htbp]
  \centering
  \begin{tabular}{l|l|l|l}
    \hline
    Kernel Name       & SL (dim)                   & DL (dim)     & Target Value (dim)                                \\
    \hline
    LapPGrad          & $q$ (1)                    & $d_j$ (3)    & $\phi,\nabla \phi$ (1+3)                          \\
    LapPGradGrad      & $q$ (1)                    & $d_j$ (3)    & $\phi,\nabla \phi, \nabla\nabla \phi$ (1+3+6).    \\
    LapQPGradGrad     & $Q_{ij}$ (9)               & --           & $\phi,\nabla \phi, \nabla\nabla \phi$ (1+3+6).    \\
    Stokeslet         & $f_j$ (3)                  & --           & $u_i$ (3)                                         \\
    RPY               & $f_j,b$ (3+1)              & --           & $u_i,\nabla^2 u_i$ (3+3)                          \\
    StokesRegVel      & $f_j,\epsilon$ (3+1)       & --           & $u_i$ (3)                                         \\
    StokesRegVelOmega & $f_k,t_l,\epsilon$ (3+3+1) & --           & $u_i,w_j$ (3+3)                                   \\
    PVel              & $f_k,q$ (3+1)              & $D_{kl}$ (9) & $p,u_i$ (1+3)                                     \\
    PVelGrad          & $f_k,q$ (3+1)              & $D_{kl}$ (9) & $p,u_i,\dpone{p}{x_i},\dpone{u_i}{x_j}$ (1+3+3+9) \\
    PVelLapLacian     & $f_k,q$ (3+1)              & $D_{kl}$ (9) & $p,u_i,\nabla^2 u_{i}$ (1+3+3)                    \\
    Traction          & $f_k,q$ (3+1)              & $D_{kl}$ (9) & $\sigma_{ij}$ (9)                                 \\
    \hline
  \end{tabular}
  \caption{
    \label{tab:implement}
    All kernels currently supported by the \texttt{STKFMM} package.
    All are supported for open boundary conditions and for singly, doubly, and triply periodic boundary conditions.
    In addition, the Stokeslet and RPY kernels are supported with a no-slip wall combined with open, singly, or doubly periodic boundary conditions.
    SL:~single-layer source points.
    DL:~double-layer source points.
  }
\end{table}

\section{Overhead of periodic boundary condition}
\label{app:pbc}
\begin{figure}[htbp]
  \includegraphics[width=\linewidth]{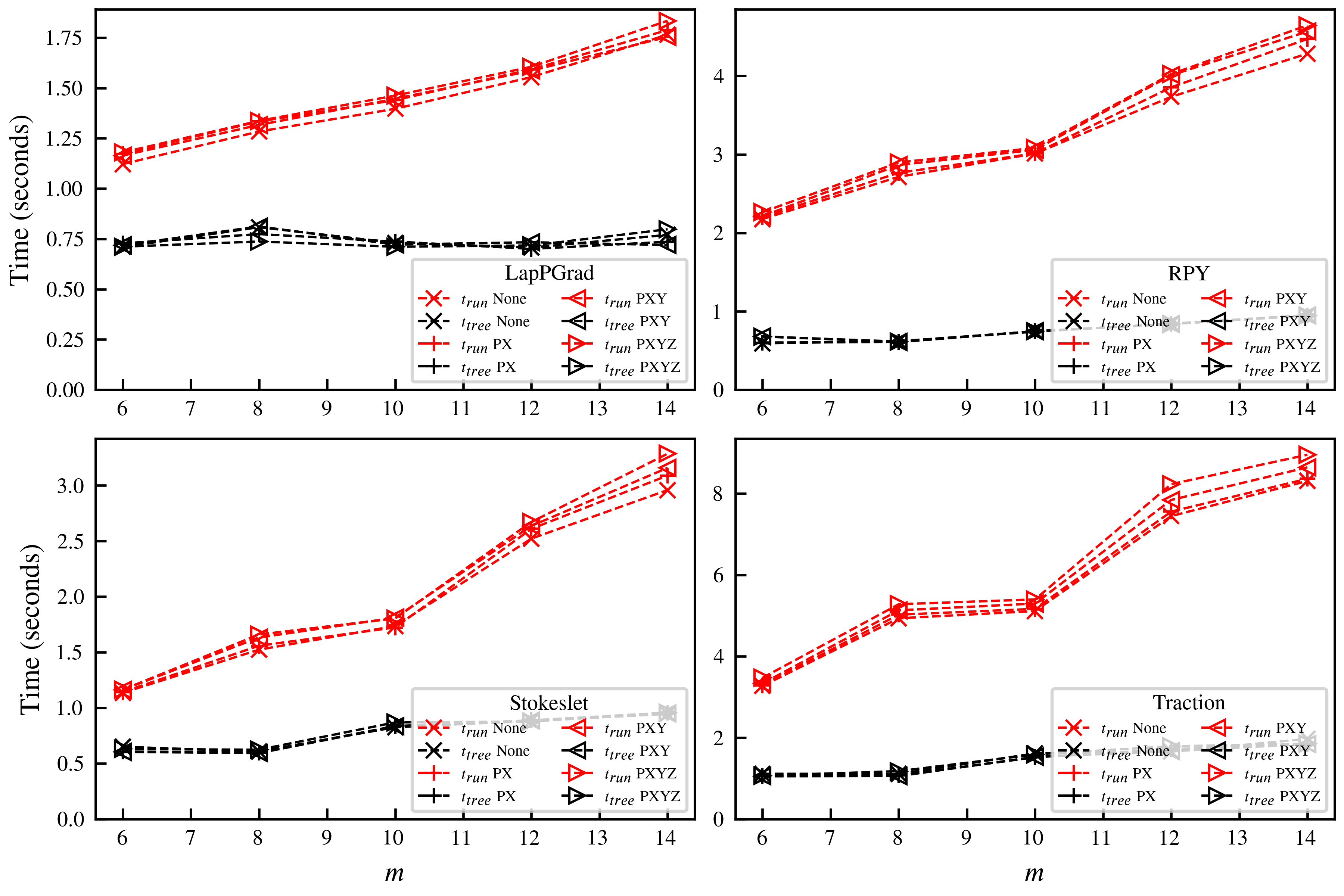}
  \caption{
    \label{fig:timepbc}
    Timing results for selected kernels with different periodic boundary conditions.
    None, PX, PXY, and PXYZ refer to open boundary, singly periodic, doubly periodic, and triply periodic, respectively.
  }
\end{figure}

In contrast to our previous implementation \cite{YanShelley2018Flexiblyimposingperiodicity}, we have merged our postfix of applying the \ML{} operator to implement periodic boundary condition into the downward pass (\LL{} and \LT{}) of running KIFMM itself, with the help of Dr. Malhotra.
Therefore, the extra overhead of running KIFMM with periodicity compared to open boundary condition is significantly reduced.
Fig.~\ref{fig:timepbc} shows the comparison of tree construciton time and FMM evaluation time for four different kernels under different boundary conditions.
The overhead of implementing periodic boundary conditions is usually $5\%$ of the evalution time with free boundary condition.

\section{Note about our SIMD implementation}
\label{app:simd}
Although the benchmarks presented in Sec.~\ref{sec:bench} were on CPUs (Intel Skylake and Cascade Lake) that support AVX-512 instructions, our kernels are hand-optimized with AVX2 instructions.
Since the AVX-512 instruction set is still evolving and new subsets of instructions are still being added, we plan to support AVX-512 in the future when the instruction set is more stable.

\section{Parallel scaling}
\label{app:scaling}
Detailed parallel scaling benchmarks and analysis are beyond the scope of this paper because we implemented our package based on the \texttt{PVFMM} library without modifying any of its advanced parallelization facilities.
Since the parallel scaling of \texttt{PVFMM} has been thoroughly tested~\cite{MalhotraBiros2015PVFMMParallelKernel}, we present only some brief test results here.
We tested strong and weak scaling for the low-precision case \(m=8\) and the high-precision case \(m=12\).
The scaling benchmarks were measured on a cluster where each node was equipped with two Intel Cascade Lake 8268 CPUs @ 2.9 GHz.
Each CPU had 24 cores, hyper-threading was disabled, and the machines were set to ``performance'' mode.
We launched two MPI ranks on each CPU.
Each MPI rank launched 12 OpenMP threads and each thread was bound to one CPU core.
The number of source and target points was increased compared to the convergence tests in the previous sections, and we increased the box size to \(L=100\).
We again used a nonuniform distribution of points picked from the log-normal distribution in Eq.~(\ref{eq:ptdist}).

\begin{figure}[htbp]
  \centering
  \includegraphics[width=\linewidth]{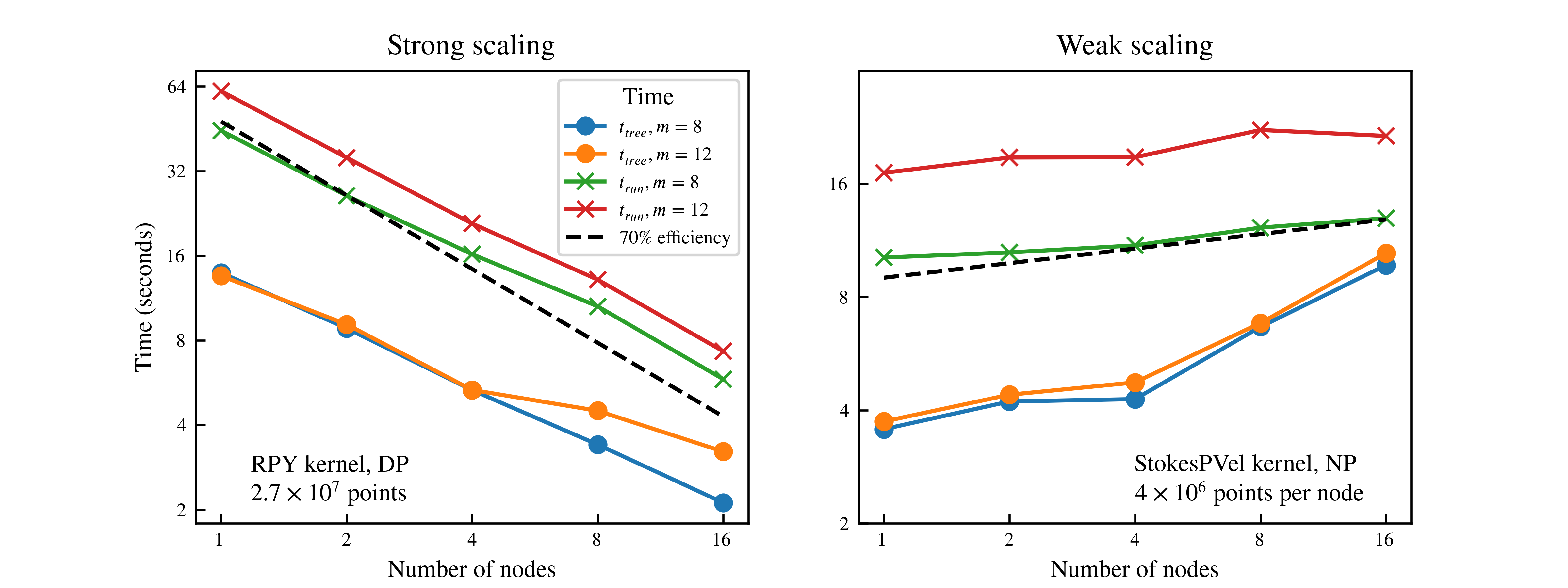}
  \caption{\label{fig:scaling}
    Strong and weak scaling benchmarks for our \texttt{KAFMM} package.
    Each node has 48 cores.
    The left panel shows strong scaling results for a fixed-size RPY tensor problem with \(2.7\times 10^7\) source and target points, with doubly periodic (DP) boundary conditions.
    The right panel shows weak scaling results for a variable-size StokesPVel problem with \(4\times 10^6\) points per node, with non-periodic (NP) boundary conditions.
    In both cases the dashed line is a reference line of \(70\%\) parallel efficiency.
  }
\end{figure}

We tested strong scaling using the RPY kernel with a doubly periodic boundary condition and \(2.7\times 10^7\) target points and source points.
Fig.~\ref{fig:scaling}(a) shows the results, which are comparable to 70\% parallel efficiency (dashed black line).
We tested weak scaling using the StokesPVel kernel \(\bG^P\) with 4 million target points per node, i.e., 1 million points for each MPI rank. As shown in Fig.~\ref{fig:scaling}(b), the \(t_\text{run}\) results are again comparable to 70\% parallel efficiency.
These results are completely determined by the \texttt{PVFMM} library, and are close to its benchmarks \cite{MalhotraBiros2015PVFMMParallelKernel}.

\bibliographystyle{apsrev4-2}
% \bibliography{ref}%

%

\end{document}